\documentclass[twoside,12pt]{article}
\usepackage{amssymb}
\usepackage{amsthm}
\usepackage{amsmath}
\usepackage[all]{xy}
\usepackage{enumerate}
\usepackage{fancyhdr}
\usepackage{a4wide}

\newtheorem{theorem}{Theorem}[section]
\newtheorem{proposition}[theorem]{Proposition}
\newtheorem{corollary}[theorem]{Corollary}
\newtheorem{lemma}[theorem]{Lemma}
\theoremstyle{definition}
\newtheorem*{notation}{Notation}
\newtheorem*{Beweis}{Proof}
\newtheorem{definition}[theorem]{Definition}
\newtheorem{punto}[theorem]{}
\theoremstyle{remark}
\newtheorem{remark}[theorem]{Remark}
\newtheorem{ex}[theorem]{Example}
\newtheorem{exs}[theorem]{Examples}
\newtheorem{remarks}[theorem]{Remarks}
\CompileMatrices
\swapnumbers
\CompileMatrices
\input{tcilatex}
\begin{document}

\title{Exact Sequences in Non-Exact Categories\\
(An Application to Semimodules)}
\author{\textbf{Jawad Y. Abuhlail\thanks{%
The author would like to acknowledge the support provided by the Deanship of
Scientific Research (DSR) at King Fahd University of Petroleum $\&$ Minerals
(KFUPM) for funding this work through project No. FT100004.}} \\
Department of Mathematics and Statistics\\
King Fahd University of Petroleum $\&$ Minerals\\
abuhlail@kfupm.edu.sa}
\date{}
\maketitle

\begin{abstract}
We consider a notion of \emph{exact sequences} in any -- not necessarily
exact -- pointed category relative to a given $(\mathbf{E},\mathbf{M})$%
-factorization structure. We apply this notion to introduce and investigate
a new notion of exact sequences of semimodules over semirings relative to
the canonical image factorization. Several homological results are proved
using the new notion of exactness including some restricted versions of the
Short Five Lemma and the Snake Lemma opening the door for introducing and
investigating \emph{homology objects} in such categories. Our results apply
in particular to the variety of commutative monoids extending results in
homological varieties to relative homological varieties.
\end{abstract}

\section*{Introduction}

\qquad \emph{Exact sequences} and \emph{exact functors} are important tools
in Homological Algebra which was developed first in the categories of
modules over rings \cite{CE1956} and generalized later to arbitrary Abelian
categories (e.g. \cite{Hel1958}). Different sets of axioms characterizing
\emph{additive} abstract categories which can be considered -- in some sense
-- \emph{natural home} for exact sequences were developed over time; such
categories were called \emph{exact} (e.g. \emph{Buchsbaum-exact categories}
\cite{Buc1955}, \emph{Quillen-exact categories} \cite{Qui1973}). For these
categories, the defining axioms are usually based on a distinguished class
of sequences, called an \emph{exact structure}, which is used to define the
(short and long) exact sequences in the resulting exact category as well as
exact functors between such exact categories. On the other hand, the
so-called \emph{Barr-exact categories }\cite{Bar1971}, which are \emph{%
regular categories }with canonical $(\mathbf{RegEpi},\mathbf{Mono})$%
-factorization structures (e.g. \cite{Gri1971}, \cite[14.E]{AHS2004}, \cite%
{Bor1994b}), provide an alternative notion of exactness in possibly \emph{%
non-additive} categories. In such categories, the role of exact sequences is
played by the so-called \emph{exact forks} which are also used to define
exact functors between Barr-exact categories. For a systematic study and
comprehensive exposition of these and other notions of exact categories, the
interested reader is advised to consult \cite{Bue2010}.

An elegant notion of exact categories to which we refer often in this
manuscript is duo to Puppe \cite{Pup1962} (see also Mitchell \cite{Mit1965}%
). We call a category $\mathfrak{C}$ a \emph{Puppe-exact category}\textit{\ }%
iff it is pointed (i.e. $\mathrm{Hom}_{\mathfrak{C}}(A,B)$ has a zero
morphism for each $A,B\in \mathrm{Obj}(\mathfrak{C})$) and has a $(\mathbf{%
NormalEpi},\mathbf{NormalMono})$\emph{-factorization structure} (e.g. \cite[%
14.F]{AHS2004}); such a category is \emph{additive} if and only if it is
Abelian (cf. \cite[3.2]{BP1969}). By \cite[13.1.3]{Sch1972}, any Puppe-exact
category has kernels and cokernels; moreover it is \emph{normal} (i.e. every
monomorphism is a kernel) and \emph{conormal} (i.e. every epimorphism is a
cokernel). The \emph{image} (\emph{coimage}) of a morphism $\gamma $ in a
Puppe-exact category $\mathfrak{C}$ is defined as $\func{Im}(\gamma ):=%
\mathrm{Ker}(\mathrm{coker}(\gamma )$ ($\mathrm{Coim}(\gamma ):=\mathrm{Coker%
}(\mathrm{ker}(\gamma ))$) and a sequence $A\overset{f}{\longrightarrow }B%
\overset{g}{\longrightarrow }C$ in $\mathfrak{C}$ is said to be \emph{exact}
iff $\func{Im}(f)\simeq \mathrm{Ker}(g)$ or equivalently $\mathrm{Coim}%
(g)\simeq \mathrm{Coker}(f)$ \cite[12.4.9, 13.1.3]{Sch1972}.

Many interesting pointed categories are not Puppe-exact, (e.g. some
varieties of Universal Algebra like the variety $\mathbf{Grp}$ of groups,
the variety $\mathbf{Mon}$ of monoids and the variety $\mathbf{pSet}$ of
pointed sets). Thus, the following question arises naturally:%
\begin{equation*}
\text{\textbf{Question:}}\mathbf{\ }\emph{When\ is\ an\ exact\ sequence\ }A%
\overset{f}{\longrightarrow }B\overset{g}{\longrightarrow }C\emph{\ in\ a\
pointed\ category\ exact?}
\end{equation*}%
The main goal of this article is providing an answer to the above mentioned
question. Our approach is based on analyzing the notion of exact sequences
in Puppe-exact categories and then generalizing it to any pointed category $%
\mathfrak{C}$ relative to a given $(\mathbf{E},\mathbf{M})$\emph{%
-factorization structure}, which always exists \cite[Section 14]{AHS2004}
(see also \cite{Bar2002}): we say that a sequence $A\overset{f}{%
\longrightarrow }B\overset{g}{\longrightarrow }C$ is $(\mathbf{E},\mathbf{M}%
) $\emph{-exact} iff there exist $f^{\prime }\in \mathbf{E}$ and $g^{\prime
\prime }\in \mathbf{M}$ such that $f=\mathrm{ker}(g)\circ f^{\prime }$ and $%
g=g^{\prime \prime }\circ \mathrm{coker}(f)$ are the \emph{essentially unique%
} $(\mathbf{E},\mathbf{M})$-factorizations of $f$ and $g$ in $\mathfrak{C}.$
To illustrate this notion of exactness, we introduce a restricted version of
the Short Exact Lemma and introduce a class of \emph{relative homological
categories} which generalizes the notion of \emph{homological categories} in
the sense of Borceux and Bourn \cite[Chapter 4]{BB2004}.

Before we proceed, we find it suitable to include the following
clarification. A successful notion of exact sequences already exists in
several pointed categories which are not Puppe-exact (e.g. in $\mathbf{Grp}$%
):\ A sequence $A\overset{f}{\longrightarrow }B\overset{g}{\longrightarrow }%
C $ of groups is exact iff $\func{Im}(f)\simeq \mathrm{Ker}(g).$ While used
in many papers, and even considered \emph{standard}, this notion of
exactness is not necessarily appropriate in other pointed categories (e.g.
in $\mathbf{Mon}$). We briefly demonstrate why we believe this is the case.
Firstly, one should be careful about the definition of the image (coimage)
of a morphisms $\gamma $ in a category which is not Puppe-exact: although
several authors define $\func{Im}(\gamma ):=\mathrm{Ker}(\mathrm{coker}%
(\gamma ))$ ($\mathrm{Coim}(\gamma ):=\mathrm{Coker}(\mathrm{ker}(\gamma ))$%
), this might not be the appropriate notion in a category which is not
Puppe-exact as it does not necessarily satisfy the universal property that
an image (coimage) is supposed to satisfy (cf. \cite[5.8.7]{Fai1973} and
\cite{EW1987}). Secondly, even if the appropriate image (coimage) is used,
one has to take into consideration the natural dual condition of exactness,
namely $\mathrm{Coim}(g)\simeq \mathrm{Coker}(f).$ This \emph{hidden}
condition is equivalent to $\func{Im}(f)\simeq \mathrm{Ker}(g)$ in
Puppe-exact categories \cite[Lemma 13.1.4]{Sch1972}; however, this is not
necessarily the case in categories which are not Puppe-exact. So, one might
end up with two different notions: \emph{left-exact sequences} for which $%
\func{Im}(f)\simeq \mathrm{Ker}(g)$ and \emph{right-exact sequences} for
which $\mathrm{Coim}(g)\simeq \mathrm{Coker}(f),$ while exact sequences have
to be defined as those which are left-exact and right-exact.

An example that demonstrates how adopting the definition of exact sequences
in Puppe-exact categories to arbitrary pointed categories might create
serious problems is the notion of exact sequences of semimodules over
semirings due to Takahashi \cite{Tak1981}. An unfortunate choice of a notion
of exactness and an inappropriate choice of a \emph{tensor functor} which is
not left adjoint of the Hom functor, in addition to the \textit{bad} nature
of monoids (in contrast with the \textit{good} nature of groups), are among
the main reasons for failing to develop a satisfactory homological theory
for semimodules or commutative monoids so far (there are indeed many
successful investigations related to the homology of monoids, e.g. \cite%
{KKM2000}).

This manuscript is divided as follows. After this introduction, and for the
convention of the reader, we recall in Section 1 some terminology and
notions from Category Theory. In particular, we analyze the notion of exact
sequences in \emph{Puppe-exact categories} and use that analysis to
introduce a new notion of exact sequences in arbitrary pointed categories.
Moreover, we present some special classes of morphisms which play an
important role in the sequel. In Section 2, we collect some definitions and
results on semirings and semimodules and clarify the differences between the
terminology used in this paper and the classical terminology; we also
clarify the reason for changing some terminology. In Section 3, we apply our
general definition of exactness to obtain a new notion of exact sequences of
semimodules over semirings. We demonstrate how this notion enables us to
characterize in a very simple way, similar to that in homological
categories, different classes of morphisms (e.g. monomorphisms, regular
epimorphisms, isomorphisms). In Section 4, we illustrate the advantage of
our notion of exactness over the existing ones by showing how it enables us
to prove some of the elementary diagram lemmas for semimodules over
semirings. Moreover, we introduce a restricted version of the \emph{Short
Five Lemma \ref{short-5}}, which characterizes the homological categories
among the pointed regular ones, and use it to introduce a new class of \emph{%
relative homological categories} w.r.t. a given factorization structure and
a special class of morphisms. The category of cancellative semimodules over
semirings, in particular the category of cancellative commutative monoids,
is introduced as a prototype of such categories. Moreover, we prove a
restricted version of the \emph{Snake Lemma \ref{snake}} for cancellative
semimodules (cancellative commutative monoids) which opens the door for
introducing and investigating \emph{homology objects} in such categories.

\section{Exact Sequences in Pointed Categories}

\qquad Throughout, and unless otherwise explicitly mentioned, $\mathfrak{C}$
is an arbitrary \emph{pointed} category (i.e. $\mathrm{Hom}_{\mathfrak{C}%
}(A,B)$ has a zero morphism); all objects and morphisms are assumed to be in
$\mathfrak{C}.$ When clear from the context, we may drop $\mathfrak{C}.$ Our
main references in Category Theory are \cite{AHS2004} and \cite{Mac1998}.

\begin{punto}
A \emph{monomorphism} in $\mathfrak{C}$ is a morphism $m$ such that for any
morphisms $f_{1},f_{2}:$%
\begin{equation*}
m\circ f_{1}=m\circ f_{2}\Rightarrow f_{1}=f_{2}.
\end{equation*}%
An \emph{equalizer} of a family of morphisms $(f_{\lambda }:A\rightarrow
B)_{\Lambda }$ in $\mathfrak{C}$ is a morphism $g:A^{\prime }\rightarrow A$
in $\mathfrak{C}$ such that $f_{\lambda }\circ g=f_{\lambda ^{\prime }}\circ
g$ for all $\lambda ,\lambda ^{\prime }\in \Lambda $ and whenever there
exists $g^{\prime }:A^{\prime \prime }\rightarrow A$ with $f_{\lambda }\circ
g^{\prime }=f_{\lambda ^{\prime }}\circ g^{\prime }$ for all $\lambda
,\lambda ^{\prime }\in \Lambda $ then there exists a \emph{unique} morphism $%
\widetilde{g}:A^{\prime \prime }\rightarrow A^{\prime }$ such that $g\circ
\widetilde{g}=g^{\prime }:$
\begin{equation*}
\xymatrix{& A'' \ar@{.>}_{\tilde{g}}[dl] \ar^{g'}[d] & & \\ A' \ar_{g}[r] &
A \ar@<1ex>^{f_{\lambda}}[r] \ar[r] & B}
\end{equation*}%
With $\mathrm{Equ}((f_{\lambda })_{\lambda \in \Lambda })$ we denote the
domain of the \emph{essentially unique} equalizer of $(f_{\lambda
})_{\lambda \in \Lambda },$ if it exists. A morphism $g$ in in $\mathfrak{C}$
is said to be a \emph{regular monomorphism} iff $g=\mathrm{equ}(f_{1},f_{2})$
for two morphisms $f_{1},f_{2}$ in $\mathfrak{C}.$
\end{punto}

\begin{punto}
Let $g$ be a morphism in $\mathfrak{C}.$ We call $\mathrm{ker}(f):=\mathrm{%
Equ}(f,0)$ the \emph{kernel} of $f.$ We say that $g$ is a \emph{normal
monomorphism} iff $g=\mathrm{ker}(f)$ for some morphism $f$ in $\mathfrak{C}%
. $ The category $\mathfrak{C}$ is said to be \emph{normal} iff every
monomorphism in $\mathfrak{C}$ is normal.
\end{punto}

\begin{punto}
An \emph{epimorphism} in $\mathfrak{C}$ is a morphism $e$ such that for any
morphisms $f_{1},f_{2}:$%
\begin{equation*}
f_{1}\circ e=f_{2}\circ e\Rightarrow f_{1}=f_{2}.
\end{equation*}%
A coequalizer of a family of morphisms $(f_{\lambda }:A\rightarrow
B)_{\Lambda }$ in $\mathfrak{C}$ is a morphism $g:B\rightarrow B^{\prime }$
in $\mathfrak{C}$ such that $g\circ f_{\lambda }=g\circ f_{\lambda ^{\prime
}}$ for all $\lambda ,\lambda ^{\prime }\in \Lambda $ and whenever there
exists $g^{\prime }:B\rightarrow B^{\prime \prime }$ with $g^{\prime }\circ
f_{\lambda }=g^{\prime }\circ f_{\lambda ^{\prime }}$ for all $\lambda
,\lambda ^{\prime }\in \Lambda $ then there exists a \emph{unique} morphism $%
\widetilde{g}:B^{\prime }\rightarrow B^{\prime \prime }$ such that $%
\widetilde{g}\circ g=g^{\prime }:$%
\begin{equation*}
\xymatrix{ A \ar@<1ex>^{f_{\lambda}}[r] \ar[r] & B \ar_{g'}[d] \ar^{g}[r] &
B' \ar@{.>}^{\tilde{g}}[dl] \\& B'' & &}
\end{equation*}%
With $\mathrm{Coequ}((f_{\lambda })_{\lambda \in \Lambda })$ we denote the
codomain of the essentially unique coequalizer of $(f_{\lambda })_{\lambda
\in \Lambda },$ if it exists. A morphism $g$ is said to be a \emph{regular
epimorphism} iff $g=\mathrm{Coequ}(f_{1},f_{2})$ for two morphisms $%
f_{1},f_{2}$ in $\mathfrak{C}.$
\end{punto}

\begin{punto}
We call $\mathrm{Coker}(f):=\mathrm{Coequ}(f,0)$ the \emph{cokernel} of $f.$
A morphism $g$ is said to be a \emph{conormal epimorphism} iff $g=\mathrm{%
coker}(f)$ for some morphism $f$ in $\mathfrak{C}.$ The category $\mathfrak{C%
}$ is said to be \emph{conormal} iff every epimorphism in $\mathfrak{C}$ is
conormal.
\end{punto}

\begin{notation}
We fix some notation:

\begin{itemize}
\item With $\mathbf{Mono}$\textbf{$($}$\mathfrak{C)}$ ($\mathbf{RegMono}$%
\textbf{$($}$\mathfrak{C)}$) we denote the class of (regular) monomorphisms
in $\mathfrak{C}$ and by $\mathbf{Epi}$\textbf{$($}$\mathfrak{C)}$ ($\mathbf{%
RegEpi}$\textbf{$($}$\mathfrak{C)}$) the class of (regular) epimorphisms in $%
\mathfrak{C}.$ We denote by $\mathbf{NormMono}$\textbf{$($}$\mathfrak{C)}%
\subseteq \mathbf{RegMono(\mathfrak{C)}\ }$($\mathbf{NormEpi(\mathfrak{C)}}%
\subseteq \mathbf{RegEpi(\mathfrak{C)}}$) the class of normal monomorphisms
(normal epimorphisms) in $\mathfrak{C}.$

\item With $\mathbf{Iso}$\textbf{$($}$\mathfrak{C)}$ we denote the class of
isomorphisms and with $\mathbf{Bimor}$\textbf{$($}$\mathfrak{C)}$ the class
of bimorphisms (i.e. monomorphisms and epimorphisms) in $\mathfrak{C}.$

\item Let $\mathfrak{C}$ be concrete (over the category $\mathbf{Set}$ of
sets) with underlying functor $U:\mathfrak{C}\longrightarrow \mathbf{Set}.$
We denote by $\mathbf{Inj}$\textbf{$($}$\mathfrak{C)}$ ($\mathbf{Surj(}%
\mathfrak{C)}$) the class of morphisms $\gamma $ in $\mathfrak{C}$ such that
$U(\gamma )$ is an injective (surjective)\ map.
\end{itemize}
\end{notation}

\begin{remark}
(e.g. \cite[7.76]{AHS2004}) We have%
\begin{equation*}
\mathbf{Iso(\mathfrak{C)}}\subseteq \mathbf{NormMono(}\mathfrak{C)}\subseteq
\mathbf{RegMono(\mathfrak{C)}}\subseteq \mathbf{Mono(\mathfrak{C)}}
\end{equation*}%
and%
\begin{equation*}
\mathbf{Iso(\mathfrak{C}})\subseteq \mathbf{NormEpi(}\mathfrak{C)}\subseteq
\mathbf{RegEpi(\mathfrak{C)}}\subseteq \mathbf{Epi(\mathfrak{C)}}.
\end{equation*}
\end{remark}

\begin{definition}
(Compare with \cite[5.8.7]{Fai1973}, \cite{EW1987}, \cite{Bar2002}) Let $%
\mathbf{E}$ and $\mathbf{M}$ be classes of morphisms in $\mathbf{\mathfrak{C}%
}$ and $\gamma :X\longrightarrow Y$ a morphism in $\mathfrak{C}.$

\begin{enumerate}
\item The $\mathbf{M}$\emph{-image} of $\gamma $ is $\mathrm{im}(\gamma ):%
\func{Im}(\gamma )\longrightarrow Y$ in $\mathbf{M}$ such that $\gamma =%
\mathrm{im}(\gamma )\circ \iota _{\gamma }$ for some morphism $\iota
_{\gamma }$ and if $\gamma =m\circ \iota $ for some $m\in \mathbf{M},$ then
there exists a unique morphism $\alpha _{m}:\func{Im}(\gamma
)\longrightarrow Z$ such that $m\circ \alpha _{m}=\mathrm{im}(\gamma )$ and $%
\alpha _{m}\circ \iota _{\gamma }=\iota .$

\item The $\mathbf{E}$\emph{-coimage} of $\gamma $ is $\mathrm{coim}(\gamma
):X\longrightarrow \mathrm{Coim}(\gamma )$ in $\mathbf{E}$ such that $\gamma
=c_{\gamma }\circ \mathrm{coim}(\gamma )$ for some morphism $c_{\gamma }$
and if $\gamma =c\circ e$ for some $e\in \mathbf{E},$ then there exists a
unique morphism $\beta _{e}:Z\longrightarrow \mathrm{Coim}(\gamma )$ such
that $\beta _{e}\circ e=\mathrm{coim}(\gamma )$ and $c_{\gamma }\circ \beta
_{e}=c.$%
\begin{equation*}
\begin{tabular}{lll}
$\xymatrix{ & & {\rm Im}(\gamma) \ar@{.>}^{{\rm im}(\gamma)}[ddddrr]
\ar@{-->}_(.65){\alpha_m}[dd] & & \\ & & & & \\ & & Z \ar_{m}[ddrr] & & \\
\\ X \ar_{\iota}[uurr] \ar_(.6){\gamma}[rrrr] \ar^{\iota_{\gamma}}[uuuurr] &
& & & Y}$ &  & $\xymatrix{ & & Z \ar@{-->}_(.65){\beta_e}[dd]
\ar^{c}[ddddrr] & & \\ & & & & \\ & & {\rm Coim}(\gamma)
\ar_{c_{\gamma}}[ddrr] & & \\ \\ X \ar^{e}[uuuurr] \ar_(.6){\gamma}[rrrr]
\ar@{.>}_{{\rm coim}(\gamma)}[uurr] & & & & Y}$%
\end{tabular}%
\end{equation*}
\end{enumerate}
\end{definition}

\begin{definition}
(\cite[14.1]{AHS2004})\ Let $\mathbf{E}$ and $\mathbf{M}$ be classes of
morphisms in $\mathbf{\mathfrak{C}}.$ The pair $(\mathbf{E},\mathbf{M})$ is
called a \emph{factorization structure }(\emph{for morphisms} in) $\mathbf{%
\mathfrak{C}},$ and $\mathbf{\mathfrak{C}}$ is said to be $(\mathbf{E},%
\mathbf{M})$\emph{-structured,} provided that

\begin{enumerate}
\item $\mathbf{E}$ and $\mathbf{M}$ are closed under composition with
isomorphisms.

\item $\mathbf{\mathfrak{C}}$ has $(\mathbf{E},\mathbf{M})$-factorizations,
\emph{i.e.} each morphism $f$ in $\mathbf{\mathfrak{C}}$ has a factorization
$f=m\circ e$ with $m\in \mathbf{M}$ and $e\in \mathbf{E}.$

\item $\mathbf{\mathfrak{C}}$ has the \emph{unique} $(\mathbf{E},\mathbf{M})$%
\emph{-diagonalization property} (or the \emph{diagonal-fill-in property})
\emph{i.e. }for each commutative square%
\begin{equation*}
\xymatrix{A \ar^{e}[r] \ar_{f}[d] & B \ar^{g}[d] \ar@{.>}^ {d}[ld] \\ C
\ar_{m}[r] & D }
\end{equation*}%
with $e\in \mathbf{E}$ and $m\in \mathbf{M},$ there exists a \emph{unique}
morphism $d:B\longrightarrow C$ such that $d\circ e=f$ and $m\circ d=g.$
\end{enumerate}
\end{definition}

\begin{punto}
Let $\mathfrak{C}$ be an $(\mathbf{E},\mathbf{M})$-structured category and $%
\gamma :X\longrightarrow Y$ be a morphism in $\mathfrak{C}$ with $(\mathbf{E}%
,\mathbf{M})$-factorization $\gamma :X\overset{e}{\longrightarrow }U\overset{%
m}{\longrightarrow }Y.$ Let $\mathrm{Coim}(\gamma )$ and $\func{Im}(\gamma )$
be the the $\mathbf{E}$-coimage of $\gamma $ and the $\mathbf{M}$-image of $%
\gamma ,$ respectively. Then there exist isomorphisms $\mathrm{Coim}(\gamma )%
\overset{d_{1}}{\simeq }U\overset{d_{2}}{\simeq }\func{Im}(\gamma )$ such
that $d_{2}\circ d_{1}$ is the canonical morphism $d_{\gamma }:\mathrm{Coim}%
(\gamma )\longrightarrow \func{Im}(\gamma ),$ which is in this case an
isomorphism:%
\begin{equation*}
\xymatrix{ & & {\rm Coim}(\gamma) \ar@{-->}_(.65){d_{\gamma}}[dd]
\ar^{c_{\gamma}}[ddddrr] & & \\ & & & & \\ & & {\rm Im}(\gamma)
\ar@{.>}_{{\rm im}(\gamma)}[ddrr] & & \\ \\ X \ar@{.>}^(.6){{\rm
coim}(\gamma)}[uuuurr] \ar_{\iota_{\gamma}}[uurr] \ar_{\gamma}[rrrr] & & & &
Y}
\end{equation*}
\end{punto}

\begin{remarks}
\begin{enumerate}
\item \label{unique}For any category, $(\mathbf{Iso},\mathbf{Mor})$ and $(%
\mathbf{Mor},\mathbf{Iso})$ are \emph{trivial} factorization structures.

\item Some authors assume that $\mathbf{E}\subseteq \mathbf{Epi}(\mathfrak{C}%
)$ and $\mathbf{M}\subseteq \mathbf{Mon}(\mathfrak{C})$ (e.g. \cite{Bar2002}%
).

\item If $(\mathbf{E},\mathbf{M})$ is a factorization structure for $%
\mathfrak{C},$ then $\mathbf{E}\cap \mathbf{M}=\mathbf{Iso}(\mathfrak{C}).$

\item As a result of the unique diagonalization property, any $(\mathbf{E},%
\mathbf{M})$-factorization in an $(\mathbf{E},\mathbf{M})$-structured
category is \emph{essentially} unique (compare with \cite[Proposition 14.4]%
{AHS2004}). Suppose that $m_{1}\circ e_{1}=\gamma =m_{2}\circ e_{2}$ are two
$(\mathbf{E},\mathbf{M})$-factorizations of a morphism $\gamma
:A\longrightarrow B$ in $\mathfrak{C}$%
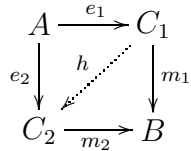
\begin{equation*}
\xymatrix{A \ar^{e_1}[r] \ar_{e_2}[d] & C_1 \ar^{m_1}[d] \ar@{.>}_{h}[dl]\\
C_2 \ar_{m_2}[r] & B}
\end{equation*}%
Then there exists a (unique) isomorphism $h:C_{1}\longrightarrow C_{2}$ s.t.
the above diagram commutes.
\end{enumerate}
\end{remarks}

\subsection*{Exact Categories}

\qquad There are several notions of \emph{exact sequences} and \emph{exact
categories} in the literature (e.g. \cite{Buc1955}, \cite{Qui1973}, \cite%
{Pup1962}, \cite{Bar1971}).

\begin{punto}
Call $\mathfrak{C}$ a\emph{\ Puppe-exact category} iff it is pointed and has
a $(\mathbf{NormalEpi},\mathbf{NormalMono})$-factorization structure. By
\cite[14.F (a)]{AHS2004}, a pointed category is Puppe-exact if and only if
it has $(\mathbf{NormalEpi},\mathbf{NormalMono})$\emph{-factorizations},
i.e. every morphism $\gamma $ admits a -- necessarily \emph{unique} --
factorization $\gamma =\gamma ^{\prime \prime }\circ \gamma ^{\prime }$ such
that $\gamma ^{\prime }$ is a cokernel and $\gamma ^{\prime \prime }$ is a
kernel. The image and the coimage of a morphism $\gamma :X\longrightarrow Y$
in such a categories are given by $\func{Im}(\gamma ):=\mathrm{Ker}(\mathrm{%
coker}(\gamma ))$ and $\mathrm{Coim}(\gamma ):=\mathrm{Coker}(\mathrm{ker}%
(\gamma )),$ respectively. Moreover, a sequence $A\overset{f}{%
\longrightarrow }B\overset{g}{\longrightarrow }C$ is said to be \emph{exact}
iff $\func{Im}(f)\simeq \mathrm{Ker}(g).$
\end{punto}

\begin{remarks}
\begin{enumerate}
\item Every non-empty Puppe-exact category has a zero-object, kernels and
cokernels, is normal, conormal and has equalizers.

\item Let $\mathfrak{C}$ be a category with a zero-object, kernels,
cokernels and equalizers. If $\mathfrak{C}$ is normal, then $\mathrm{Coker}(%
\mathrm{ker}(\gamma ))\simeq \mathrm{Ker}(\mathrm{coker}(\gamma ))$ for any
morphism $\gamma $ in $\mathfrak{C}$ \cite[Proposition 5.20]{Fai1973}.
\end{enumerate}
\end{remarks}

\qquad In light of the previous remarks, \cite[Lemma 31.14]{Sch1972} can be
restated as follows:

\begin{lemma}
\label{P-exact}Let $\mathfrak{C}$ be a Puppe-exact category, $A\overset{f}{%
\longrightarrow }B\overset{g}{\longrightarrow }C$ a sequence in $\mathfrak{C}
$ with $g\circ f=0$ and consider the following commutative diagram with
canonical and induced factorizations
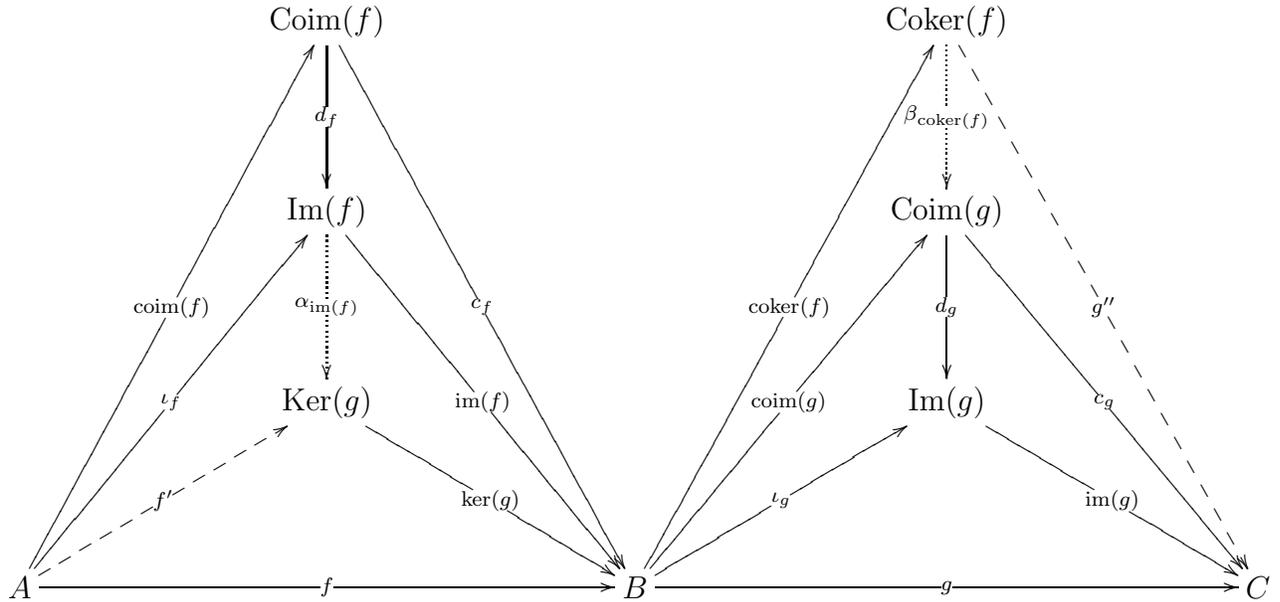
\begin{equation*}
\xymatrix{ & & & {\rm Coim}(f) \ar|-{d_f}[dd] \ar|-{c_f}[ddddddrrr] & & & &
& & {\rm Coker}(f) \ar@{.>}|-{\beta_{{\rm coker}(f)}}[dd]
\ar@{-->}|-{g''}[ddddddrrr] & & & \\ & & & & & & & & & & & & \\ & & & {\rm
Im}(f) \ar@{.>}|-{\alpha_{{\rm im}(f)}}[dd] \ar|-{{\rm im}(f)} [ddddrrr] & &
& & & & {\rm Coim}(g) \ar|-{d_g}[dd] \ar|-{c_g}[ddddrrr] & & & \\ & & & & &
& & & & & & & \\ & & & {\rm Ker}(g) \ar[ddrrr]|-{{\rm ker}(g)} & & & & & &
{\rm Im}(g) \ar[ddrrr]|-{{\rm im}(g)} & & & \\ & & & & & & & & & & & & \\ A
\ar@{-->}|-{f'}[uurrr] \ar|-{f}[rrrrrr] \ar|-{\iota_f}[rrruuuu]
\ar[uuuuuurrr]|-{{\rm coim}(f)} & & & & & & B \ar|-{g}[rrrrrr]
\ar[uuuuuurrr]|-{{\rm coker}(f)} \ar[rrruuuu]|-{{\rm coim}(g)}
\ar|-{\iota_g}[uurrr] & & & & & & C }
\end{equation*}%
The following are equivalent:

\begin{enumerate}
\item $A\overset{f}{\longrightarrow }B\overset{g}{\longrightarrow }C$ is
exact \emph{(}i.e. $\func{Im}(f)\overset{\alpha _{\mathrm{im}(f)}}{\simeq }%
\mathrm{Ker}(g)$\emph{)};

\item $\mathrm{Coker}(f)\overset{\beta _{\mathrm{coker}(f)}}{\simeq }\mathrm{%
Coim}(g);$

\item $\mathrm{Coim}(f)\simeq \mathrm{Ker}(g);$

\item $\mathrm{Coker}(f)\simeq \func{Im}(g);$

\item $\mathrm{Im}(f)\simeq \mathrm{Ker}(\mathrm{coim}(g));$

\item $\mathrm{Coim}(g)\simeq \mathrm{Coker}(\mathrm{im}(f));$

\item $f^{\prime }$ is a (normal) epimorphism;

\item $g^{\prime \prime }$ is a (normal) monomorphism.
\end{enumerate}
\end{lemma}

Inspired by the previous lemma, we introduce a notion of exact sequences in
any pointed category:

\begin{definition}
Let $\mathfrak{C}$ be any pointed category and fix a factorization structure
$(\mathbf{E},\mathbf{M})\ $for $\mathfrak{C}.$ We call a sequence%
\begin{equation}
A\overset{f}{\longrightarrow }B\overset{g}{\longrightarrow }C  \label{ABC}
\end{equation}%
\emph{exact} w.r.t. $(\mathbf{E},\mathbf{M})$ iff $f$ and $g$ have
factorizations
\begin{equation*}
f=\mathrm{ker}(g)\circ f^{\prime }\text{ and }g=g^{\prime \prime }\circ
\mathrm{coker}(f)\text{ with }(f^{\prime },\mathrm{ker}(g)),(\mathrm{coker}%
(f),g^{\prime \prime })\in \mathbf{E}\times \mathbf{M}.
\end{equation*}%
\begin{equation*}
\xymatrix{& A \ar@{.>}_{f'}[dl] \ar^{f}[d] & \\ {\rm Ker}(g) \ar_{{\rm
ker}(g)}[r] & B \ar_{g}[d] \ar^(0.4){{\rm coker}(f)}[r] & {\rm Coker}(f)
\ar@{.>}^{g''}[dl] \\ & C & }
\end{equation*}%
When the factorization structure is clear from the context we drop it. We
call a sequence%
\begin{equation*}
\cdots \longrightarrow A_{i-1}\overset{f_{i-1}}{\longrightarrow }A_{i}%
\overset{f_{i}}{\longrightarrow }A_{i+1}\longrightarrow \cdots
\end{equation*}%
\emph{exact} at $A_{i}$ iff $A_{i-1}\overset{f_{i-1}}{\longrightarrow }A_{i}%
\overset{f_{i}}{\longrightarrow }A_{i+1}$ is exact; moreover, we call this
sequence \emph{exact} iff it is exact at $A_{i}$ for every $i.$ An exact
sequence%
\begin{equation}
0\longrightarrow A\overset{f}{\longrightarrow }B\overset{g}{\longrightarrow }%
C\longrightarrow 0  \label{ses}
\end{equation}%
is called a \emph{short exact sequence}.
\end{definition}

\begin{remark}
Let $\mathfrak{C}$ be a pointed category and fix a factorization structure $(%
\mathbf{E},\mathbf{M})\ $for $\mathfrak{C}.$ It follows immediately from the
definition that (\ref{ses}) is a short exact sequence w.r.t. $(\mathbf{E},%
\mathbf{M})$ if and only if $\mathrm{Coker}(f)\in \mathbf{E},$ $\mathrm{Ker}%
(g)\in \mathbf{M},$ $f=\mathrm{Ker}(g)$ and $g=\mathrm{Coker}(f).$
\end{remark}

\begin{ex}
Let $\mathfrak{C}$ be a Puppe exact category. A sequence $A\overset{f}{%
\longrightarrow }B\overset{g}{\longrightarrow }C$ is \emph{exact} if and
only if $f=\mathrm{ker}(g)\circ f^{\prime },$ $g=g^{\prime \prime }\circ
\mathrm{coker}(f)$ with $(f^{\prime },\mathrm{ker}(g)),$ $(\mathrm{coker}%
(f),g^{\prime \prime })\in \mathbf{NormalEpi}\times \mathbf{NormalMono}.$
Notice that, by Lemma \ref{P-exact}, this is equivalent to the classical
notion of exacts sequences in Puppe-exact categories, namely $\func{Im}%
(f)\simeq \mathrm{Ker}(g).$ This applies in particular to the categories of
modules over rings (e.g. the category $\mathbf{Ab}$ of Abelian groups).
\end{ex}

\begin{ex}
Let $\mathfrak{C}$ be a pointed $(\mathbf{E},\mathbf{M})$-structure category
with $\mathbf{NormalEpi}\subseteq \mathbf{E}\subseteq \mathbf{Epi}$ and $%
\mathbf{NormalMono}\subseteq \mathbf{M}\subseteq \mathbf{Mono}.$ Then a
sequence $A\overset{f}{\longrightarrow }B\overset{g}{\longrightarrow }C$ in $%
\mathfrak{C}$ is exact if and only if the essentially unique $(\mathbf{E},%
\mathbf{M})$ factorizations $f=m_{1}\circ e_{1},$ $g=m_{2}\circ e_{2}$ can
be chosen so that $m_{1}=\mathrm{ker}(e_{2})$ and $e_{2}=\mathrm{coker}%
(m_{1}).$ Moreover, a sequence $0\longrightarrow A\overset{f}{%
\longrightarrow }B\overset{g}{\longrightarrow }C\longrightarrow 0$ is exact
if and only if $f=\mathrm{Ker}(g)$ and $g=\mathrm{Coker}(f).$ This applies
to general pointed categories which are $(\mathbf{RegEpi},\mathbf{Mono})$%
-structured or $(\mathbf{Epi},\mathbf{RegMono})$-structured. In particular,
this applies to pointed regular categories (compare with \cite[Definition
4.1.7]{BB2004}).
\end{ex}

\begin{ex}
Let $\mathfrak{C}$ be a pointed protomodular category (in the sense of D.
Bourn \cite{Bou1991}) with finite limits. By \cite[Proposition 3.1.23]%
{BB2004}, $g\in \mathbf{RegEpi}(\mathfrak{C})$ if and only if $g=\mathrm{%
coker}(\mathrm{ker}(g)).$ If $\mathfrak{C}$ is $(\mathbf{RegEpi},\mathbf{Mono%
})$-structured or $(\mathbf{Epi},\mathbf{RegMono})$-structured, then it
follows that a sequence $0\longrightarrow A\overset{f}{\longrightarrow }B%
\overset{g}{\longrightarrow }C\longrightarrow 0$ in $\mathfrak{C}$ is exact
if and only if $f=\mathrm{ker}(g)$ and $g$ is a regular epimorphism. This
applies in particular to homological categories, which are precisely pointed
and protomodular regular categories \cite{BB2004}.
\end{ex}

\begin{ex}
Let $(\mathfrak{C};\mathbf{E})$ be a \emph{relative homological category} in
the sense of \cite{Jan2006}, where $\mathbf{E}$ is a distinguished class of
normal epimorphisms and assume that $\mathfrak{C}$ is $(\mathbf{E},\mathbf{%
Mono})$-structured (which is not actually assumed in the defining axioms of
such categories). Analyzing Condition (a) on $g^{\prime }$ (page 192), which
was assumed to prove the so called \emph{Relative Snake Lemma, }shows that
this assumption and along with the assumptions on $f^{\prime }$ are
essentially equivalent to assuming that the row $0\longrightarrow A^{\prime }%
\overset{f^{\prime }}{\longrightarrow }B^{\prime }\overset{g^{\prime }}{%
\longrightarrow }C^{\prime }$ is $(\mathbf{E},\mathbf{Mono})$-exact.
\end{ex}

\subsection*{Steady Morphisms}

\qquad In what follows, we consider a special class of categories to which
there is a natural transfer of the notion of exact sequences in Puppe-exact
categories.

\begin{definition}
Let $\mathfrak{C}$ be a pointed $(\mathbf{E},\mathbf{M})$-structured
category. We say that a morphism $\gamma :X\longrightarrow Y$ in $\mathfrak{C%
}$ is:

\emph{steady} w.r.t. $(\mathbf{E},\mathbf{M})$ iff $\mathrm{Ker}(\gamma ),$ $%
\mathrm{Coker}(\mathrm{ker}(\gamma ))$ exist in $\mathfrak{C}$ and $\gamma $
admits an $(\mathbf{E},\mathbf{M})$- factorization $\gamma =\gamma ^{\prime
\prime }\circ \mathrm{coker}(\mathrm{ker}(\gamma )),$ equivalently $\mathrm{%
Coker}(\mathrm{ker}(\gamma ))\simeq \mathrm{Coim}(\gamma );$

\emph{costeady} w.r.t. $(\mathbf{E},\mathbf{M})$ iff $\mathrm{Coker}(\gamma
),$ $\mathrm{Ker}(\mathrm{coker}(\gamma ))$ exist in $\mathfrak{C}$ and $%
\gamma $ admits an $(\mathbf{E},\mathbf{M})$- factorization $\gamma =\mathrm{%
ker}(\mathrm{coker}(\gamma ))\circ \gamma ^{\prime },$ equivalently $\mathrm{%
Ker}(\mathrm{coker}(\gamma ))\simeq \func{Im}(\gamma );$

\emph{bisteady} w.r.t. $(\mathbf{E},\mathbf{M})$ iff $\gamma $ is steady and
costeady w.r.t. $(\mathbf{E},\mathbf{M}),$ equivalently $\mathrm{Coker}(%
\mathrm{ker}(\gamma ))\simeq \mathrm{Coim}(\gamma )\overset{d_{\gamma }}{%
\simeq }\func{Im}(\gamma )\simeq \mathrm{Ker}(\mathrm{ker}(\gamma )).$
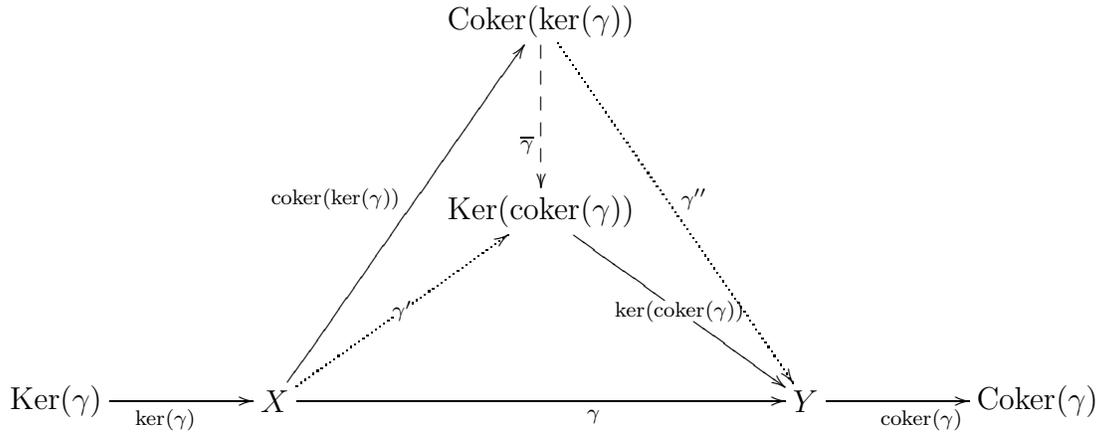
\begin{equation*}
\xymatrix{ & & & & {\rm Coker}({\rm ker}(\gamma)) \ar@{.>}^{\gamma
''}[ddddrr] \ar@{-->}_(.65){\overline{\gamma}}[dd] & & & & \\ & & & & & & &
& \\ & & & & {\rm Ker} ({\rm coker}(\gamma)) \ar|-{{\rm ker} ({\rm
coker}(\gamma))}[ddrr] & & & & \\ & & & & & & & & \\ {\rm Ker}(\gamma)
\ar_{{\rm ker}(\gamma)}[rr] & & X \ar@{.>}|-{\gamma '}[uurr]
\ar_(.6){\gamma}[rrrr] \ar^{{\rm coker} ({\rm ker}(\gamma))}[uuuurr] & & & &
Y \ar_{{\rm coker}(\gamma)}[rr] & & {\rm Coker}(\gamma)}
\end{equation*}%
We call $\mathfrak{C}$ \emph{steady }(resp. \emph{costeady, bisteady})
w.r.t. $(\mathbf{E},\mathbf{M})$ iff all morphisms in $\mathfrak{C}$ are
steady (resp. costeady, bisteady) w.r.t. $(\mathbf{E},\mathbf{M}).$
\end{definition}

\begin{remark}
Let $\mathfrak{C}$ be a pointed $(\mathbf{E},\mathbf{M})$-structured
category. If $\mathfrak{C}$ is bisteady w.r.t. $(\mathbf{E},\mathbf{M}),$
then $\mathfrak{C}$ is Puppe-exact: in this case, every morphism in $%
\mathfrak{C}$ has a $(\mathbf{NormalEpi},\mathbf{NormalMono})$%
-factorization, whence $\mathfrak{C}$ is Puppe-exact \cite[14.F]{AHS2004}.
Moreover, if $\mathbf{NormalEpi}\subseteq \mathbf{E}$ and $\mathbf{NormalMono%
}\subseteq \mathbf{M}$ then $\mathfrak{C}$ is bisteady w.r.t. $(\mathbf{E},%
\mathbf{M})$ if and only if $\mathfrak{C}$ is Puppe-exact.
\end{remark}

\begin{punto}
\label{UA}All varieties -- in the sense of Universal Algebra -- are $(%
\mathbf{RegEpi},\mathbf{Mono})$-structured. Moreover, the class of regular
epimorphisms coincides with that of surjective morphisms, and the class of
monomorphisms coincides with that of injective morphisms. Let $\mathcal{V}$
be a pointed variety. We say that a morphism $\gamma :X\longrightarrow Y$ in
$\mathcal{V}$ is \emph{steady} (resp. \emph{costeady}, \emph{bisteady}) iff $%
\gamma $ is steady (resp. costeady, bisteady) w.r.t. $(\mathbf{Surj},\mathbf{%
Inj}).$ With $\func{Im}(\gamma )$ ($\mathrm{Coim}(\gamma )$) we mean the $%
\mathbf{Inj}$-image (the $\mathbf{Surj}$-coimage) of $\gamma .$ Moreover, we
say that a sequence $X\overset{f}{\longrightarrow }Y\overset{g}{%
\longrightarrow }Z$ in $\mathcal{V}$ is \emph{exact} iff it is $(\mathbf{Surj%
},\mathbf{Inj})$-exact.
\end{punto}

\begin{ex}
The variety $\mathbf{Grp}$ of all (Abelian and non-Abelian) groups is
steady. Let $\gamma :X\longrightarrow Y$ be any morphism of groups. Notice
that $\mathrm{Ker}(\gamma )=\{x\in X\mid \gamma (x)=1_{Y}\}$ while $\mathrm{%
Coker}(\gamma )=Y/N_{\gamma },$ where $N_{\gamma }$ is the \emph{normal
closure} of $\gamma (X).$ Consider the canonical $(\mathbf{Surj},\mathbf{Inj}%
)$-factorization $\gamma :X\overset{\mathrm{im}(\gamma )}{\longrightarrow }%
\gamma (X)\overset{\iota }{\longrightarrow }Y$ where $\iota $ is the
canonical embedding. Consider also the factorization $\gamma :X\overset{%
\mathrm{coker}(\mathrm{ker}(\gamma ))}{\longrightarrow }X/\mathrm{Ker}%
(\gamma )\overset{\gamma ^{\prime \prime }}{\longrightarrow }Y.$ Clearly, $%
\gamma $ is steady if and only if $\gamma ^{\prime \prime }$ is injective.
Indeed, if $\gamma ^{\prime \prime }([x_{1}])=\gamma ^{\prime \prime
}([x_{2}])$ for some $x_{1},x_{2}\in X,$ then $\gamma (x_{1})=\gamma (x_{2})$
whence $\gamma (x_{1}^{-1}x_{2})=1_{Y}$ and it follows that $%
x_{1}^{-1}x_{2}=k$ for some $k\in \mathrm{Ker}(\gamma ),$ i.e. $%
[x_{1}]=[x_{2}].$ Consequently, $\gamma ^{\prime \prime }\in \mathbf{Inj}.$
On the other hand, consider the factorization $\gamma :X\overset{\gamma
^{\prime }}{\longrightarrow }N_{\gamma }\overset{\mathrm{ker}(\mathrm{coker}%
(\gamma ))}{\longrightarrow }Y.$ Then $\gamma $ is costeady if and only if $%
\gamma (X)=N_{\gamma }$ if and only if $\gamma (X)\leq G$ is a normal
subgroup. Clearly, $\mathbf{Grp}$ is not costeady: Let $G$ be a group, $H$ a
subgroup that is not normal in $G$ and let $\gamma :H\hookrightarrow G$ be
the embedding. Indeed, $H=\gamma (H)\neq N_{\gamma },$ \emph{i.e.} $\gamma $
is not costeady. Consequently, $\mathbf{Grp}$ is not a bisteady category.
\end{ex}

\qquad In the following example, we demonstrate how the \emph{classical }%
notion of exact sequences of groups is consistent with our new definition of
exact sequences in arbitrary pointed categories.

\begin{ex}
\label{G-exact}Let $A\overset{f}{\longrightarrow }B\overset{g}{%
\longrightarrow }C$ be a sequence of groups and consider the canonical
factorizations of $f:A\overset{\mathrm{im}(f)}{\longrightarrow }f(A)\overset{%
\iota _{1}}{\hookrightarrow }B$ and $g:B\overset{\mathrm{im}(g)}{%
\longrightarrow }g(B)\overset{\iota _{2}}{\hookrightarrow }C.$ If the given
sequence is exact, then $f=\mathrm{ker}(g)\circ f^{\prime }$ with $f^{\prime
}$ surjective. This implies that $f(A)=\mathrm{Ker}(g).$ On the other hand,
assume that $f(A)=\mathrm{Ker}(g).$ Then $f$ has a an $(\mathbf{Inj},\mathbf{%
Surj})$-factorization as $f=\mathrm{ker}(g)\circ \mathrm{im}(f).$ Moreover,
it is evident that there is an isomorphism of groups $B/\mathrm{Ker}(g)%
\overset{\gamma }{\simeq }g(B).$ So, $g$ has an $(\mathbf{Inj},\mathbf{Surj}%
) $-factorization $g=(\iota _{2}\circ \gamma )\circ \mathrm{coker}(g).$ It
follows that $A\overset{f}{\longrightarrow }B\overset{g}{\longrightarrow }C$
is exact if and only if $f(A)=\mathrm{Ker}(g).$
\end{ex}

\section{Semirings and Semimodules}

\qquad \emph{Semirings} (\emph{semimodules}) are roughly speaking, rings
(modules) without subtraction. Semirings were studied independently by
several algebraists, especially by H. S. Vandiver \cite{Van1934} who
considered them as the \emph{best} algebraic structures which unify rings
and bounded distributive lattices. Since the sixties of the last century,
semirings were shown to have significant applications in several areas as
Automata Theory (e.g. \cite{Eil1974}, \cite{Eil1976}, \cite{KS1986}),
Theoretical Computer Science (e.g. \cite{HW1998}) and many classical areas
of mathematics (e.g. \cite{Go19l99a}, \cite{Gol1999b}).

Recently, semirings played an important role in several emerging areas of
research like Idempotent Analysis (e.g. \cite{KM1997}, \cite{LMS2001}, \cite%
{Lit2007}), Tropical Geometry (e.g. \cite{R-GST2005}, \cite{Mik2006}) and
many aspects of modern Mathematics and Mathematical Physics (e.g. \cite%
{Gol2003}, \cite{LM2005}). In his dissertation \cite{Dur2007}, N. Durov
demonstrated that semirings are in one-to-one correspondence with what he
called \emph{algebraic additive monads} on the category $\mathbf{Set}$ of
sets. Moreover, a connection between semirings and the so-called $\mathbb{F}$%
-rings, where $\mathbb{F}$ is the field with one element, was pointed out in
\cite[1.3 -- 1.4]{PL}.

The theory of semimodules was developed mainly by M. Takahashi, who
published several fundamental papers on this topic (cf. \cite{Tak1979} --
\cite{Tak1985}) and to whom research in the theory of semimodules over
semirings is indeed indebted. However, it seems to us that there are some
gaps in his theory of semimodules which led to confusion and sometimes
conceptual misunderstandings. Instead of introducing appropriate definitions
and notions that fit well with the special properties of the category of
semimodules over semirings, some definitions and notions which are fine in
\emph{Puppe-exact categories }in general, and in categories of modules over
rings in particular, were enforced in a category which is, in general, far
away from being Puppe-exact.

A systematic development of the homological theory of semirings and
semimodules has been initiated recently in a series of papers by Y. Katsov
\cite{Kat1997} and carried out in a continuing series of papers (e.g. \cite%
{Kat2004a}, \cite{Kat2004b}, \cite{KTN2009}, \cite{KN2011}, \cite{IK2011},
\cite{IK2011}). Another approach that is worth mentioning was initiated by
A. Patchkoria in \cite{Pat1998} and continued in a series of papers (e.g.
\cite{Pat2000a}, \cite{Pat2000b}, \cite{Pat2003}, \cite{Pat2006}, \cite%
{Pat2009}).

In what follows, we revisit the category of semimodules over a semiring. In
particular, we adopt a new definition of exact sequences of semimodules and
investigate it. We also introduce some terminology that will be needed in
the sequel.

\begin{punto}
\label{semig}Let $(S,\ast )$ be a semigroup. We call $s\in S$ \emph{%
cancellable} iff for any $s_{1},s_{2}\in S:$%
\begin{equation*}
s_{1}\ast s=s_{2}\ast s\Longrightarrow s_{1}=s_{2}\text{ and }s\ast
s_{1}=s\ast s_{2}\Longrightarrow s_{1}=s_{2}.
\end{equation*}%
We call $S$ \emph{cancellative} iff all elements of $S$ are cancellable. We
say that a morphism of semigroups $f:S\longrightarrow S^{\prime }$ is \emph{%
cancellative} iff $f(s)\in S^{\prime }$ is cancellable for every $s\in S.$
We call $S$ an \emph{idempotent semigroup} iff $s\ast s=s$ for every $s\in
S. $
\end{punto}

\begin{punto}
Let $(S,+)$ be an Abelian additive semigroup. A subset $X\subseteq S$ is
said to be \emph{subtractive }iff for any $s\in S$ and $x\in X$ we have: $%
s+x\in X\Longrightarrow s\in X.$ The\emph{\ subtractive closure} of a
non-empty subset $X\subseteq S$ is given by%
\begin{equation*}
\overline{X}:=\{s\in S\mid s+x_{1}=x_{2}\text{ for some }x_{1},x_{2}\in X\}.
\end{equation*}%
If $X$ is a subsemigroup of $S,$ then indeed $X$ is subtractive if and only
if $X=\overline{X}.$ We call a morphism of Abelian semigroups $%
f:S\longrightarrow S^{\prime }$ \emph{subtractive} iff $f(S)\subseteq
S^{\prime }$ is subtractive, equivalently iff%
\begin{equation*}
f(S)=\{s^{\prime }\in S^{\prime }\mid \text{ }s^{\prime }+f(s_{1})=f(s_{2})%
\text{ for some }s_{1},s_{2}\in S\}.
\end{equation*}
\end{punto}

\begin{punto}
A \emph{semiring} is an algebraic structure $(S,+,\cdot ,0,1)$ consisting of
a non-empty set $S$ with two binary operations \textquotedblleft $+$%
\textquotedblright\ (addition) and \textquotedblleft $\cdot $%
\textquotedblright\ (multiplication) satisfying the following conditions:

\begin{enumerate}
\item $(S,+,0)$ is an Abelian monoid with neutral element $0_{S};$

\item $(S,\cdot ,1)$ is a monoid with neutral element $1;$

\item $x\cdot (y+z)=x\cdot y+x\cdot z$ and $(y+z)\cdot x=y\cdot x+z\cdot x$
for all $x,y,z\in S;$

\item $0\cdot s=0=s\cdot 0$ for every $s\in S$ (i.e. $0$ is \emph{absorbing}%
).
\end{enumerate}

Let $S,S^{\prime }$ be semirings. A map $f:S\rightarrow S^{\prime }$ is said
to be a \emph{morphism of semirings} iff for all $s_{1},s_{2}\in S:$%
\begin{equation*}
f(s_{1}+s_{2})=f(s_{1})+f(s_{2}),\text{ }f(s_{1}s_{2})=f(s_{1})f(s_{2}),%
\text{ }f(0_{S})=0_{S^{\prime }}\text{ and }f(1_{S})=1_{S^{\prime }}.
\end{equation*}%
The category of semirings is denoted by $\mathbf{SRng}.$
\end{punto}

\begin{punto}
Let $(S,+,\cdot )$ be a semiring. We say that $S$ is

\emph{cancellative} iff the additive semigroup $(S,+)$ is cancellative;

\emph{commutative} iff the multiplicative semigroup $(S,\cdot )$ is
commutative;

\emph{semifield} iff $(S\backslash \{0\},\cdot ,1)$ is a commutative group.
\end{punto}

\begin{exs}
Rings are indeed semirings. A trivial, but important, example of a \emph{%
commutative }semiring is $(\mathbb{N}_{0},+,\cdot )$ (the set of
non-negative integers). Indeed, $(\mathbb{R}_{0}^{+},+,\cdot )$ and $(%
\mathbb{Q}_{0}^{+},+,\cdot )$ are semifields. A more interesting example is
the semi-ring $(\mathrm{ideal}(R),+,\cdot )$ consisting of all ideals of a
(not necessarily commutative) ring; this appeared first in the work of \emph{%
Dedekind} \cite{Ded1894}. On the other hand, for an integral domain $R,$ $(%
\mathrm{ideal}(R),+,\cap )$ is a semiring if and only if $R$ is a Pr\"{u}fer
domain. Every bounded distributive lattice $(R,\vee ,\wedge )$ is a
commutative (additively) idempotent semiring. The \emph{additively idempotent%
} semirings $\mathbb{R}_{\max }:=(\mathbb{R}\cup \{-\infty \},\max ,+)$ and $%
\mathbb{R}_{\min }:=(\mathbb{R}\cup \{\infty \},\min ,+)$ play an important
role in idempotent and tropical mathematics (e.g. \cite{Lit2007}); the
subsemirings $\mathbb{N}_{\max }:=(\mathbb{N}\cup \{-\infty \},\max ,+)$ and
$\mathbb{N}_{\min }:=(\mathbb{N}\cup \{\infty \},\min ,+)$ played an
important role in Automata Theory (e.g. \cite{Eil1974}, \cite{Eil1976}). The
singleton set $S=\{0\}$ is a semiring with the obvious addition and
multiplication. In the sequel, we always assume that $0_{S}\neq 1_{S}$ so
that $S\neq \{0\},$ the \emph{zero semiring}.
\end{exs}

\begin{punto}
Let $S$ be a semiring. A \emph{right }$S$\emph{-semimodule} is an algebraic
structure $(M,+,0_{M};\leftharpoondown )$ consisting of a non-empty set $M,$
a binary operation \textquotedblleft $+$\textquotedblright\ along with a
right $S$-action%
\begin{equation*}
M\times S\longrightarrow M,\text{ }(m,s)\mapsto ms,
\end{equation*}%
such that:

\begin{enumerate}
\item $(M,+,0_{M})$ is an Abelian monoid with neutral element $0_{M};$

\item $(ms)s^{\prime }=m(ss^{\prime }),$ $(m+m^{\prime })s=ms+m^{\prime }s$
and $m(s+s^{\prime })=ms+ms^{\prime }$ for all $s,s^{\prime }\in S$ and $%
m,m^{\prime }\in M;$

\item $m1_{S}=m$ and $m0_{S}=0_{M}=0_{M}s$ for all $m\in M$ and $s\in S.$

Let $M,M^{\prime }$ be right $S$-semimodules. A map $f:M\rightarrow
M^{\prime }$ is said to be a \emph{morphism of right }$S$\emph{-semimodules}
(or $S$\emph{-linear}) iff for all $m_{1},m_{2}\in M$ and $s\in S:$%
\begin{equation*}
f(m_{1}+m_{2})=f(m_{1})+f(m_{2})\text{ and }f(ms)=f(m)s.
\end{equation*}%
The set $\mathrm{Hom}_{S}(M,M^{\prime })$ of $S$-linear maps from $M$ to $%
M^{\prime }$ is clearly a monoid under addition. The category of right $S$%
-semimodules is denoted by $\mathbb{S}_{S}.$ Similarly, one can define the
category of left $S$-semimodules $_{S}\mathbb{S}.$ A right $S$-semimodule $%
M_{S}$ is said to be \emph{cancellative} iff the semigroup $(M,+)$ is
cancellative. With $\mathbb{CS}_{S}\subseteq \mathbb{S}_{S}$ (resp. $_{S}%
\mathbb{CS}\subseteq $ $_{S}\mathbb{S})$ we denote the full subcategory of
cancellative right (left) $S$-semimodules.
\end{enumerate}
\end{punto}

\begin{punto}
Let $M$ be a right $S$-semimodule. A non-empty subset $L\subseteq M$ is said
to be an $S$\emph{-subsemimodule, }and we write $L\leq _{S}M,$ iff $L$ is
closed under \textquotedblleft $+_{M}$\textquotedblright\ and $ls\in L$ for
all $l\in L$ and $s\in S.$
\end{punto}

\begin{ex}
Every Abelian monoid $(M,+,0_{M})$ is an $\mathbb{N}_{0}$-semimodule in the
obvious way. Moreover, the categories $\mathbf{CMon}$ of commutative monoids
and the category $\mathbb{S}_{\mathbb{N}_{0}}$ of $\mathbb{N}_{0}$%
-semimodules are isomorphic.
\end{ex}

\subsection*{Congruences}

\begin{punto}
\label{zuka}Let $M$ be an $S$-semimodule. An equivalence relation
\textquotedblleft $\equiv $\textquotedblright\ on $M$ is a said to be an $S$%
\emph{-congruence on }$M$ iff for any $m,m^{\prime },m_{1},m_{1}^{\prime
},m_{2},m_{2}^{\prime }\in M$ and $s\in S$ we have
\begin{equation*}
\lbrack m_{1}\equiv m_{1}^{\prime }\text{ and }m_{2}\equiv m_{2}^{\prime
}\Rightarrow \lbrack m_{1}+m_{2}\equiv m_{1}^{\prime }+m_{2}^{\prime }]\text{
and }[m\equiv m^{\prime }\Rightarrow ms\equiv m^{\prime }s].
\end{equation*}%
The set $M/\equiv $ of equivalence classes inherit a structure of an $S$%
-semimodule in the obvious way and there is a canonical surjection of $S$%
-semimodules $\pi _{\equiv }:M\rightarrow M/\equiv .$
\end{punto}

\begin{punto}
Let $M$ be an $S$-semimodule. Every $S$-subsemimodule $L\leq _{S}M$ induces
two $S$-congruences on $M:$ the \emph{Bourne relation}%
\begin{equation*}
m_{1}\equiv _{L}m_{2}\Leftrightarrow m_{1}+l_{1}=m_{2}+l_{2}\text{ for some }%
l_{1},l_{2}\in L;
\end{equation*}%
and the \emph{Iizuka relation}%
\begin{equation*}
m_{1}[\equiv ]_{L}m_{2}\Leftrightarrow m_{1}+l_{1}+m^{\prime
}=m_{2}+l_{2}+m^{\prime }\text{ for some }l_{1},l_{2}\in L\text{ and }%
m^{\prime }\in M.
\end{equation*}%
We call the $S$-semimodule $M/L:=M/_{\equiv _{L}}$ the \emph{quotient of }$M$%
\emph{\ by }$L$ or the \emph{factor semimodule} of $M$ by $L.$ One can
easily check that $M/L=M/\overline{L}.$ If $M$ is cancellative, then $L$ and
$M/L$ are cancellative. On the other hand, the $S$-semimodule $M/[\equiv
]_{L}$ is cancellative.
\end{punto}

\begin{proposition}
The category $\mathbb{S}_{S}$ and its full subcategory $\mathbb{CS}_{S}$
have kernels and cokernels, where for any morphism of $S$-semimodules $%
f:M\rightarrow N$ we have%
\begin{equation*}
\mathrm{Ker}(f)=\{m\in M\mid f(m)=0\}\text{ and }\mathrm{\mathrm{Coker}}%
(f)=N/f(M).
\end{equation*}
\end{proposition}

Taking into account the fact that $\mathbb{S}_{S}$ is a variety (in the
sense of Universal Algebra) we have:

\begin{proposition}
\label{cc}\emph{(\cite{Tak1982b}, \cite{Tak1982c}, \cite{TW1989})}\ The
category of semimodules is

\begin{enumerate}
\item complete (i.e. has equalizers $\&$ products);

\item cocomplete (i.e. has coequalizers $\&\ $coproducts);

\item Barr-exact categories \emph{\cite{Bar1971}}.
\end{enumerate}
\end{proposition}

\begin{remark}
In \cite{Tak1982c}, Takahashi proved that the category of semimodules over a
semiring is $c$\emph{-cocomplete}, which is a \emph{relaxed} notion of
cocompleteness which he introduced. However, it was pointed to the author by
F. Linton (and other colleagues from the Category List) that such a category
is indeed cocomplete in the classical sense since it is a variety.
\end{remark}

\begin{punto}
As a variety, the category of $S$-semimodules is regular; in particular, $%
\mathbb{S}_{S}$ has a $(\mathbf{RegEpi},\mathbf{Mono})$-factorization
structure. Let $\gamma :X\longrightarrow Y$ be a morphism of $S$%
-semimodules. Then $\func{Im}(\gamma )=\gamma (X)$ and $\mathrm{Coim}(\gamma
)=X/f,$ where $X/f$ is the quotient semimodule $X/\equiv _{f}$ given by $%
x\equiv _{f}x^{\prime }$ iff $f(x)=f(x^{\prime }).$ Indeed, we have a
canonical isomorphism%
\begin{equation*}
d_{\gamma }:\mathrm{Coim}(\gamma )\simeq \func{Im}(\gamma ),\text{ }%
[x]\mapsto \gamma (x).
\end{equation*}
\end{punto}

\begin{remark}
Takahashi defined the \emph{image} of a morphism $\gamma :X\rightarrow Y$ of
$S$-semimodules as $\mathrm{Ker}(\mathrm{coker}(\gamma ))$ and the \emph{%
proper image} as $\gamma (X).$ In fact, $\gamma (X)$ is the \emph{image} of $%
\gamma $ in the categorical sense (e.g. \cite[5.8.7]{Fai1973}, \cite{EW1987}%
).
\end{remark}

\begin{punto}
We call a morphism of $S$-semimodules $\gamma :M\longrightarrow N:$

\emph{subtractive} iff $\gamma (M)\subseteq N$ is subtractive;

\emph{strong} iff $\gamma (M)\subseteq N$ is strong;

$k$\emph{-uniform} iff for any $x_{1},x_{2}\in X:$%
\begin{equation}
\gamma (x_{1})=\gamma (x_{2})\Longrightarrow \text{ }\exists \text{ }%
k_{1},k_{2}\in \mathrm{Ker}(\gamma )\text{ s.t. }x_{1}+k_{1}=x_{2}+k_{2};
\label{k-steady}
\end{equation}

$i$\emph{-uniform} iff $\gamma (X)=\overline{\gamma (X)}:=\{y\in Y\mid
y+\gamma (x_{1})=\gamma (x_{2})$ for some $x_{1},x_{2}\in X\};$

\emph{uniform }iff $\gamma $ is $k$-uniform and $i$-uniform;

\emph{semi-monomorphism} iff $\mathrm{Ker}(\gamma )=0;$

\emph{semi-epimorphism} iff $\overline{\gamma (X)}=Y;$

\emph{semi-isomorphism} iff $\mathrm{Ker}(\gamma )=0$ and$\ \overline{\gamma
(X)}=Y.$
\end{punto}

\begin{remark}
The uniform ($k$-uniform, $i$-uniform) morphisms of semimodules were called
\emph{regular} ($k$\emph{-regular, }$i$-\emph{regular}) by Takahashi \cite%
{Tak1982c}. We think that our terminology avoids confusion sine a regular
monomorphism (regular epimorphism) has a different well-established meaning
in the language of Category Theory.
\end{remark}

\begin{lemma}
\label{co-stead}Let $\gamma :X\longrightarrow Y$ be a morphism of $S$%
-semimodules.

\begin{enumerate}
\item The following are equivalent:

\begin{enumerate}
\item $\gamma $ is steady;

\item $\mathrm{Coker}(\mathrm{ker}(\gamma ))\simeq \mathrm{Coim}(\gamma );$

\item $X/\mathrm{Ker}(\gamma )\simeq \gamma (X);$

\item $\gamma $ is $k$-uniform.
\end{enumerate}

\item The following are equivalent:

\begin{enumerate}
\item $\gamma $ is costeady;

\item $\mathrm{Ker}(\mathrm{coker}(\gamma ))\simeq \func{Im}(\gamma );$

\item $\overline{\gamma (X)}=\gamma (X);$

\item $\gamma $ is $i$-uniform (subtractive).
\end{enumerate}

\item The following are equivalent:

\begin{enumerate}
\item $\gamma $ is bisteady;

\item $\mathrm{Coker}(\mathrm{ker}(\gamma ))\simeq \mathrm{Ker}(\mathrm{Coker%
}(\gamma ));$

\item $X/\mathrm{Ker}(\gamma )\simeq \overline{\gamma (X)};$

\item $\gamma $ is uniform;
\end{enumerate}
\end{enumerate}
\end{lemma}

\begin{Beweis}
Notice that the canonical $(\mathbf{Surj},\mathbf{Mono})$- factorization of $%
\gamma $ is given by $\gamma :X\overset{\mathrm{coim}(\gamma )}{%
\longrightarrow }\gamma (X)\overset{\mathrm{im}(\gamma )}{\hookrightarrow }%
Y. $

\begin{enumerate}
\item By definition, $\gamma $ is steady iff $\gamma $ admits a $(\mathbf{%
Surj},\mathbf{Mono})$- factorization $\gamma =m_{1}\circ \mathrm{coker}(%
\mathrm{ker}(\gamma )).$ It follows that $\gamma $ is steady if and only if $%
\mathrm{Coker}(\mathrm{ker}(\gamma ))\simeq \mathrm{Coim}(\gamma )$ if and
only if $X/\mathrm{Ker}(\gamma )\simeq \gamma (X)$ which is equivalent to $%
\gamma $ being $k$-uniform.

\item By definition $\gamma $ is costeady if and only if $\gamma $ admits a $%
(\mathbf{Surj},\mathbf{Mono})$- factorization $\gamma =\mathrm{ker}(\mathrm{%
coker}(\gamma ))\circ e_{2}.$ It follows that $\gamma $ is costeady if and
only if $\gamma (X)=\mathrm{Ker}(\mathrm{coker}(\gamma )).$ Notice that%
\begin{equation*}
\begin{tabular}{lll}
$\mathrm{Ker}(\mathrm{\mathrm{coker}}(\gamma ))$ & $=$ & $\{y\in Y\mid
y\equiv _{\gamma (X)}0\}$ \\
& $=$ & $\{y\in Y\mid y+\gamma (x_{1})=\gamma (x_{2})$ for some $%
x_{1},x_{2}\in X\}$ \\
& $=$ & $\overline{\gamma (X)}.$%
\end{tabular}%
\end{equation*}%
It follows that $\gamma $ is costeady if and only if $\gamma (X)=\overline{%
\gamma (X)}$ which is equivalent to $\gamma $ being subtractive.

\item This is a combination of \textquotedblleft 1\textquotedblright\ and
\textquotedblleft 2\textquotedblright .$\blacksquare $
\end{enumerate}
\end{Beweis}

\begin{punto}
Let $M$ be an $S$-semimodule, $L\leq _{S}M$ an $S$-subsemimodule and
consider the factor semimodule $M/L.$ Then we have a surjective morphism of $%
S$-semimodules%
\begin{equation*}
\pi _{L}:=M\rightarrow M/L,\text{ }m\mapsto \lbrack m]
\end{equation*}%
with%
\begin{equation*}
\mathrm{Ker}(\pi _{L})=\{m\in M\mid m+l_{1}=l_{2}\text{ for some }%
l_{1},l_{2}\in L\}=\overline{L};
\end{equation*}%
in particular, $L=\mathrm{Ker}(\pi _{L})$ if and only if $L\subseteq M$ is
subtractive.
\end{punto}

\section{Exact Sequences of Semimodules}

\qquad Throughout this section, $S$ is a ring, an $S$-semimodule is a right $%
S$-semimodule unless otherwise explicitly specified. Moreover, $\mathbb{S}%
_{S}$ ($\mathbb{CS}_{S}$) denotes the category of (cancellative) right $S$%
-semimodules.

The notion of \emph{exact sequences} of semimodules adopted by Takahashi
\cite{Tak1981} ($L\overset{f}{\longrightarrow }M\overset{g}{\longrightarrow }%
N$ is exact iff $\overline{f(M)}=\mathrm{Ker}(g)$) seems to be inspired by
the definition of exact sequences in Puppe-exact categories. We believe it
is inappropriate. The reason for this is that neither $\mathrm{Ker}(\mathrm{%
coker}(f))=\overline{f(L)}$ is the appropriate \emph{image} of $f$ nor is $%
\mathrm{Coker}(\mathrm{ker}(g))=B/\mathrm{Ker}(g)$ the appropriate \emph{%
coimage} of $g.$

Being a Barr-exact category, a natural tool to study exactness in the
category of semimodules is that of an \emph{exact fork}, introduced in \cite%
{Bar1971} and applied to study exact functors between categories of
semimodules by Katsov et al. in \cite{KN2011}. However, since the category
of semimodules has additional features, one still expects to deal with exact
sequences rather than the more complicated exact forks.

In addition to Takahashi's classical definition of exact sequences of
semimodules, two different notions of exactness for sequences of semimodules
over semirings were introduced recently. The first is due to Patchkoria \cite%
{Pat2003} ($L\overset{f}{\longrightarrow }M\overset{g}{\longrightarrow }N$
is exact iff $f(L)=\mathrm{Ker}(g)$) and the second is due to Patil and
Deore \cite{PD2006} ($L\overset{f}{\longrightarrow }M\overset{g}{%
\longrightarrow }N$ is exact iff ${\overline{f(L)}}=\mathrm{Ker}(g)$ and $g$
is \emph{steady}). Each of these definitions is stronger than Takahashi's
notion of exactness and each proved to be more efficient in establishing
some nice homological results for semimodules over semirings. However, no
clear \emph{categorical} justification for choosing either of these two
definitions was provided. A closer look at these definitions shows that they
are in fact dual to each other in some sense, and so no it not suitable --
in our opinion -- to choose one of them and drop the other. This motivated
us to introduce in Section one a new notion of exact sequences in general
pointed varieties. Applied to categories of semimodules, it turned out that
our notion of exact sequences of semimodules is in fact a combination of the
two notions of exact sequences of semimodules in the sense of \cite{Pat2003}
and \cite{PD2006}. For the sake of completeness, we mention here that there
is another notion of exact sequences of semimodules that was introduced in
\cite{AM2002}. However, the definition is rather technical and introduced
new definitions of \emph{epic} and \emph{monic} morphisms that are different
from the classical ones.

As indicated for general varieties in \ref{UA}, the category of semimodules
is $(\mathbf{RegEpi},\mathbf{Mono})$-structured, $\mathbf{RegEpi}=\mathbf{%
Surj}$ and $\mathbf{Mono}=\mathbf{Inj}.$ We say that a morphism of
semimodules $\gamma :X\longrightarrow Y$ is \emph{steady} (resp. \emph{%
costeady}, \emph{bisteady}) iff $\gamma $ is steady (resp. costeady,
bisteady) w.r.t. $(\mathbf{Surj},\mathbf{Inj}).$ Moreover, we say that a
sequence $X\overset{f}{\longrightarrow }Y\overset{g}{\longrightarrow }Z$ of
semimodules is \emph{exact} iff it is $(\mathbf{Surj},\mathbf{Inj})$-exact.

\begin{lemma}
\label{S-mon-epi}Let%
\begin{equation}
L\overset{f}{\rightarrow }M\overset{g}{\rightarrow }N  \label{lmn}
\end{equation}%
be a sequence of $S$-semimodules with $g\circ f=0$ and consider the induced
morphisms $f^{\prime }:L\rightarrow \mathrm{Ker}(g)$ and $g^{\prime \prime }:%
\mathrm{Coker}(f)\rightarrow N.$

\begin{enumerate}
\item If $f^{\prime }$ is an epimorphism, then $\overline{f(L)}=\mathrm{Ker}%
(g).$

\item $f^{\prime }$ is a regular epimorphism (surjective) if and only if $%
f(L)=\mathrm{Ker}(g)$ if and only if $\overline{f(L)}=\mathrm{Ker}(g)$ and $%
f $ is $i$-uniform.

\item $g^{\prime \prime }:\mathrm{Coker}(f)\rightarrow N$ is a monomorphism
if and only if $\overline{f(L)}=\mathrm{Ker}(g)$ and $g$ is $k$-uniform.
\end{enumerate}
\end{lemma}

\begin{Beweis}
Since $g\circ f=0,$ we have $f(L)\subseteq \overline{f(L)}\subseteq \mathrm{%
Ker}(g).$

\begin{enumerate}
\item Assume that $f^{\prime }:L\rightarrow \mathrm{Ker}(g)$ is an
epimorphism. Suppose that $\overline{f(L)}\subsetneqq \mathrm{Ker}(g),$ so
that there exist $m^{\prime }\in \mathrm{Ker}(g)\backslash f(L).$ Consider
the $S$-linear maps
\begin{equation*}
L\overset{\widetilde{f}}{\rightarrow }\mathrm{Ker}(g)\overset{f_{1}}{%
\underset{f_{2}}{\rightrightarrows }}\mathrm{Ker}(g)/f(L),
\end{equation*}%
where $f_{1}(m)=[m]$ and $f_{2}(m)=[0]$ for all $m\in \mathrm{Ker}(g).$ For
each $l\in L$ we have
\begin{equation*}
(f_{1}\circ f^{\prime })(l)=[f(l)]=[0]=(f_{2}\circ f^{\prime })(l).
\end{equation*}%
Whence, $f_{1}\circ f^{\prime }=f_{2}\circ f^{\prime }$ while $f_{1}\neq
f_{2}$ (since $f_{1}(m^{\prime })=[m^{\prime }]\neq \lbrack
0]=f_{2}(m^{\prime });$ otherwise $m^{\prime }+f(l_{1})=f(l_{2})$ for some $%
l_{1},l_{1}^{\prime }\in L$ and $m^{\prime }\in \overline{f(L)}$ which
contradicts our assumption). So, $f^{\prime }$ is not an epimorphism, a
contradiction. Consequently, $\overline{f(L)}=\mathrm{Ker}(g).$

\item Clear.

\item $(\Rightarrow )$ Assume that $g^{\prime \prime }:\mathrm{Coker}%
(f)\rightarrow N$ is a monomorphism. Let $m\in \mathrm{Ker}(g),$ so that $%
g(m)=0.$ Then $g^{\prime \prime }([m])=0.$ Since $g^{\prime \prime }$ is a
monomorphism, we have $[m]=[0]$ and so $m+f(l)=f(l^{\prime })$ for some $%
l,l^{\prime }\in L,$ whence $m\in \overline{f(L)}.$ Suppose now that $%
g(m)=g(m^{\prime })$ for some $m,m^{\prime }\in M.$ Then $g^{\prime \prime
}([m])=g^{\prime \prime }([m^{\prime }])$ and it follows, by the injectivity
of $g^{\prime \prime },$ that $[m]=[m^{\prime }]$ which implies that $%
m_{1}+m_{1}=m^{\prime }+m_{1}^{\prime }$ for some $m_{1},m_{1}^{\prime }\in
\overline{f(L)}=\mathrm{Ker}(g).$ So, $g$ is $k$-uniform.

$(\Leftarrow )$ Assume that $\overline{f(L)}=\mathrm{Ker}(g)$ and that $g$
is $k$-uniform. Suppose that $g^{\prime \prime }([m])=g^{\prime \prime
}([m^{\prime \prime }])$ for some $m_{1},m_{2}\in M.$ Then $g(m)=g(m^{\prime
}).$ Since $g$ is $k$-uniform, we $m+k=m^{\prime }+k^{\prime }$ for some $%
k,k^{\prime }\in \mathrm{Ker}(g)=\overline{f(L)}$ and it follows that $%
[m]=[m^{\prime }]$ (notice that $M/f(L)=M/\overline{f(L)}$).$\blacksquare $
\end{enumerate}
\end{Beweis}

\begin{corollary}
A sequence of semimodules $L\overset{f}{\rightarrow }M\overset{g}{%
\rightarrow }N$ is exact if and only if $f(L)=\mathrm{Ker}(g)$ and $g$ is $k$%
-uniform.
\end{corollary}

\begin{remarks}
\label{CS-mon-epi}

\begin{enumerate}
\item A morphism of cancellative semimodules $h:X\rightarrow Y$ is an
epimorphism in $\mathbb{CS}_{S}$ if and only if $\overline{h(X)}=Y.$ Indeed,
if $h$ is an epimorphism, then it follows by Lemma \ref{S-mon-epi} that $%
\overline{h(X)}=Y$ (take $g:Y\rightarrow 0$ as the zero-morphism). On the
other hand, assume that $\overline{h(L)}=Y.$ Let $Z$ be any cancellative
semimodule and consider any $S$-linear maps
\begin{equation*}
X\overset{h}{\rightarrow }Y\underset{h_{2}}{\overset{h_{1}}{%
\rightrightarrows }}Z
\end{equation*}%
with $h_{1}\circ h=h_{2}\circ h.$ Let $y\in Y$ be arbitrary. By assumption, $%
y+h(x_{1})=h(x_{2})$ for some $x_{1},x_{2}\in X,$ whence%
\begin{equation*}
h_{1}(y)+(h_{1}\circ h)(x_{1})=(h_{1}\circ h)(x_{2})=(h_{2}\circ
h)(x_{2})=h_{2}(y)+(h_{2}\circ h)(x_{1}).
\end{equation*}%
Since $Z$ is cancellative, we conclude that $h_{1}(y)=h_{2}(y).$

\item Consider the embedding $\iota :\mathbb{N}_{0}\hookrightarrow \mathbb{Z}
$ in $\mathbb{CS}_{\mathbb{N}_{0}}.$ Indeed, $\overline{\mathbb{N}_{0}}=%
\mathbb{Z},$ whence $\iota $ is an epimorphism which is not regular.

\item Let $L\overset{f}{\rightarrow }M\overset{g}{\rightarrow }N$ be a
sequence in $\mathbb{CS}_{S}$ with $g\circ f=0.$ By \textquotedblleft
1\textquotedblright , the induced morphism $f^{\prime }:L\rightarrow \mathrm{%
Ker}(g)$ is an epimorphism if and only if $\overline{f(L)}=\mathrm{Ker}(g).$
\end{enumerate}
\end{remarks}

\begin{punto}
\label{def-exact}We call a sequence of $S$-semimodules $L\overset{f}{%
\rightarrow }M\overset{g}{\rightarrow }N:$

\emph{proper-exact} iff $f(L)=\mathrm{Ker}(g);$

\emph{semi-exact} iff $\overline{f(L)}=\mathrm{Ker}(g);$

\emph{quasi-exact} iff $\overline{f(L)}=\mathrm{Ker}(g)$ and $g$ is $k$%
-uniform;

\emph{uniform} (resp. $k$\emph{-uniform}, $i$\emph{-uniform}) iff $f$ and $g$
are uniform (resp. $k$-uniform, $i$-uniform).
\end{punto}

\begin{punto}
We call a (possibly infinite) sequence of $S$-semimodules
\begin{equation}
\cdots \rightarrow M_{i-1}\overset{f_{i-1}}{\rightarrow }M_{i}\overset{f_{i}}%
{\rightarrow }M_{i+1}\overset{f_{i+1}}{\rightarrow }M_{i+2}\rightarrow \cdots
\label{chain}
\end{equation}

\emph{chain complex} iff $f_{j+1}\circ f_{j}=0$ for every $j;$

\emph{exact} (resp. \emph{proper-exact}, \emph{semi-exact}) iff each partial
sequence with three terms $M_{j}\overset{f_{j}}{\rightarrow }M_{j+1}\overset{%
f_{j+1}}{\rightarrow }M_{j+2}$ is exact (resp. proper-exact, semi-exact);

\emph{uniform }(resp. $k$-\emph{uniform}, $i$-\emph{uniform}) iff $f_{j}$ is
uniform (resp. $k$-uniform, $i$-uniform) for every $j.$
\end{punto}

\begin{definition}
Let $M$ be an $S$-semimodule.

\begin{enumerate}
\item A subsemimodule $L\leq _{S}M$ is said to be a \emph{uniform }(\emph{%
normal})\emph{\ }$S$-\emph{subsemimodule} iff the embedding $%
0\longrightarrow L\overset{\iota }{\rightarrow }M$ is uniform (normal).

\item A quotient $M/\rho ,$ where $\rho $ is an $S$-congruence relation on $%
M,$ is is said to be a \emph{uniform }(\emph{conormal})\emph{\ quotient} iff
the surijection $\pi _{L}:M\rightarrow M/\rho $ is uniform (conormal).
\end{enumerate}
\end{definition}

\begin{remark}
Every normal subsemimodule (normal quotient) is uniform.
\end{remark}

\qquad The following result can be easily verified.

\begin{lemma}
\label{i-uniform}Let $L\overset{f}{\rightarrow }M\overset{g}{\rightarrow }N$
be a sequence of semimodules.

\begin{enumerate}
\item Let $g$ be injective.

\begin{enumerate}
\item $f$ is $k$-uniform if and only if $g\circ f$ is $k$-uniform.

\item If $g\circ f$ is $i$-uniform (uniform), then $f$ is $i$-uniform
(uniform).

\item Assume that $g$ is $i$-uniform. Then $f$ is $i$-uniform (uniform) if
and only if $g\circ f$ is $i$-uniform (uniform).
\end{enumerate}

\item Let $f$ be surjective.

\begin{enumerate}
\item $g$ is $i$-uniform if and only if $g\circ f$ is $i$-uniform.

\item If $g\circ f$ is $k$-uniform (uniform), then $g$ is $k$-uniform
(uniform).

\item Assume that $f$ is $k$-uniform. Then $g$ is $k$-uniform (uniform) if
and only if $g\circ f$ is $k$-uniform (uniform).
\end{enumerate}
\end{enumerate}
\end{lemma}

\begin{Beweis}
\begin{enumerate}
\item Let $g$ be injective; in particular, $g$ is $k$-uniform.

\begin{enumerate}
\item Assume that $f$ is $k$-uniform. Suppose that $(g\circ
f)(l_{1})=(g\circ f)(l_{2})$ for some $l_{1},l_{2}\in L.$ Since $g$ is
injective, $f(l_{1})=f(l_{2}).$ By assumption, there exist $k_{1},k_{2}\in
\mathrm{Ker}(f)$ such that $l_{1}+k_{1}=l_{2}+k_{2}.$ Since $\mathrm{Ker}%
(f)\subseteq \mathrm{Ker}(g\circ f),$ we conclude that $g\circ f$ is $k$%
-uniform. On the other hand, assume that $g\circ f$ is $k$-uniform. Suppose
that $f(l_{1})=f(l_{2})$ for some $l_{1},l_{2}\in L.$ Then $(g\circ
f)(l_{1})=(g\circ f)(l_{2})$ and so there exist $k_{1},k_{2}\in \mathrm{Ker}%
(g\circ f)$ such that $l_{1}+k_{1}=l_{2}+k_{2}.$ Since $g$ is injective, $%
\mathrm{Ker}(g\circ f)=\mathrm{Ker}(f)$ whence $f$ is $k$-uniform.

\item Assume that $g\circ f$ is $i$-uniform. Let $m\in \overline{f(L)},$ so
that $m+f(l_{1})=f(l_{2})$ for some $l_{1},l_{2}\in L.$ Then $g(m)\in
\overline{(g\circ f)(L)}=(g\circ f)(L).$ Since $g$ is injective, $m\in f(L).$
So, $f$ is $i$-uniform.

\item Assume that $g$ and $f$ are $i$-uniform. Let $n\in \overline{(g\circ
f)(L)},$ so that $n+g(f(l_{1}))=g(f(l_{2}))$ for some $l_{1},l_{2}\in L.$
Since $g$ is $i$-uniform, $n\in g(M)$ say $n=g(m)$ for some $m\in M.$ But $g$
is injective, whence $m+f(l_{1})=f(l_{2}),$ i.e. $m\in \overline{f(L)}=f(L)$
since $f$ is $i$-uniform. So, $n=g(m)\in (g\circ f)(L).$ We conclude that $%
g\circ f$ is $i$-uniform.
\end{enumerate}

\item Let $f$ be surjective; in particular, $f$ is $i$-uniform.

\begin{enumerate}
\item Assume that $g$ is $i$-uniform. Let $n\in \overline{(g\circ f)(L)}$ so
that $n+g(f(l_{1}))=g(f(l_{2}))$ for some $l_{1},l_{2}\in L.$ Since $g$ is $%
i $-uniform, $n=g(m)$ for some $m\in M.$ Since $f$ is surjective, $n=g(m)\in
(g\circ f)(L).$ So, $g\circ f$ is $i$-uniform.

On the other hand, assume that $g\circ f$ is $i$-uniform. Let $n\in
\overline{g(M)},$ so that $n+g(m_{1})=g(m_{2})$ for some $m_{1},m_{2}\in M.$
Sine $f$ is surjective, there exist $l_{1},l_{2}\in L$ such that $%
f(l_{1})=m_{1}$ and $f(l_{2})=m_{2}.$ Then, $n+(g\circ f)(l_{1})=(g\circ
f)(l_{2}),$ i.e. $n\in \overline{(g\circ f)(L)}=(g\circ f)(L)\subseteq g(M).$
So, $g$ is $i$-uniform.

\item Assume that $g\circ f$ is $k$-uniform. Suppose that $g(m_{1})=g(m_{2})$
for some $m_{1},m_{2}\in M.$ Since $f$ is surjective, we have $(g\circ
f)(l_{1})=(g\circ f)(l_{2})$ for some $l_{1},l_{2}\in L.$ By assumption, $%
g\circ f$ is $k$-uniform and so there exist $k_{1},k_{2}\in \mathrm{Ker}%
(g\circ f)$ such that $l_{1}+k_{1}=l_{2}+k_{2}$ whence $%
m_{1}+f(k_{1})=m_{2}+f(k_{2}).$ Indeed, $f(k_{1}),f(k_{2})\in \mathrm{Ker}%
(g).$ i.e. $g$ is $k$-uniform.

\item Assume that $f$ and $g$ are $k$-uniform. Suppose that $(g\circ
f)(l_{1})=(g\circ f)(l_{2})$ for some $l_{1},l_{2}\in L.$ Since $g$ is $k$%
-uniform, we have $f(l_{1})+k_{1}=f(l_{2})+k_{2}$ for some $k_{1},k_{2}\in
\mathrm{Ker}(g).$ But $f$ is surjective; whence $k_{1}=f(l_{1}^{\prime })$
and $k_{2}=f(l_{2}^{\prime })$ for some $l_{1},l_{2}\in L,$ i.e. $%
f(l_{1}+l_{1}^{\prime })=f(l_{2}+l_{2}^{\prime }).$ Since $f$ is $k$%
-uniform, $l_{1}+l_{1}^{\prime }+k_{1}^{\prime }=l_{2}+l_{2}^{\prime
}+k_{2}^{\prime }$ for some $k_{1}^{\prime },k_{2}^{\prime }\in \mathrm{Ker}%
(f).$ Indeed, $l_{1}^{\prime }+k_{1}^{\prime },l_{2}^{\prime }+k_{2}^{\prime
}\in \mathrm{Ker}(g\circ f).$ We conclude that $g\circ f$ is $k$-uniform.$%
\blacksquare $
\end{enumerate}
\end{enumerate}
\end{Beweis}

\begin{remark}
Let $L\leq _{S}M\leq _{S}N$ be $S$-semimodules. It follows directly from the
previous lemma that if $L$ is uniform in $N,$ then $L$ is a uniform in $M$
as well. Moreover, if $M$ is uniform in $N,$ then $L$ is uniform in $N$ if
and only if $L$ is uniform in $M.$
\end{remark}

Our notion of exactness allows characterization of special classes of
morphisms in a way similar to that in homological categories (compare with
\cite[Proposition 4.1.9]{BB2004}, \cite[Propositions 4.4, 4.6]{Tak1981},
\cite[Proposition 15.15]{Go19l99a}):

\begin{proposition}
\label{inj-surj}Consider a sequence of semimodules%
\begin{equation*}
0\longrightarrow L\overset{f}{\longrightarrow }M\overset{g}{\longrightarrow }%
N\longrightarrow 0.
\end{equation*}

\begin{enumerate}
\item The following are equivalent:

\begin{enumerate}
\item $0\longrightarrow L\overset{f}{\rightarrow }M$ is exact;

\item $\mathrm{Ker}(f)=0$ and $f$ is steady;

\item $f$ is semi-monomorphism and $k$-uniform;

\item $f$ is injective;

\item $f$ is a monomorphism.
\end{enumerate}

\item $0\longrightarrow L\overset{f}{\longrightarrow }M\overset{g}{%
\longrightarrow }N$ is semi-exact and $f$ is uniform if and only if $L\simeq
\mathrm{Ker}(g).$

\item $0\longrightarrow L\overset{f}{\longrightarrow }M\overset{g}{%
\longrightarrow }N$ is exact if and only if $L\simeq \mathrm{Ker}(g)$ and $g$
is $k$-uniform.

\item The following are equivalent:

\begin{enumerate}
\item $M\overset{\gamma }{\rightarrow }N\rightarrow 0$ is exact;

\item $\mathrm{Coker}(\gamma )=0$ and $\gamma $ is costeady;

\item $\gamma $ is semi-epimorphism and $i$-uniform;

\item $\gamma $ is surjective;

\item $\gamma $ is a regular epimorphism;

\item $\gamma $ is a subtractive epimorphism
\end{enumerate}

\item $L\overset{f}{\rightarrow }M\overset{g}{\rightarrow }N\rightarrow 0$
is semi-exact and $g$ is uniform if and only if $N\simeq \mathrm{Coker}(f).$

\item $L\overset{f}{\longrightarrow }M\overset{g}{\longrightarrow }%
N\longrightarrow 0$ is exact if and only if $N\simeq \mathrm{Coker}(f)$ and $%
f$ is $i$-uniform.
\end{enumerate}
\end{proposition}

\begin{corollary}
The following are equivalent:

\begin{enumerate}
\item $0\rightarrow L\overset{f}{\rightarrow }M\overset{g}{\rightarrow }%
N\rightarrow 0$ is a exact sequence of $S$-semimodules;

\item $L\simeq \mathrm{Ker}(g)$ and $\mathrm{\mathrm{\mathrm{\mathrm{Coker}}}%
}(f)\simeq N;$

\item $f$ is injective, $f(L)=\mathrm{Ker}(g),$ $g$ is surjective and ($k$%
-)uniform.

In this case, $f$ and $g$ are uniform morphisms.
\end{enumerate}
\end{corollary}

\begin{remark}
A morphism of semimodules $\gamma :X\longrightarrow Y$ is an isomorphism if
and only if $0\longrightarrow X\longrightarrow Y\longrightarrow 0$ is exact
if and only if $\gamma $ is a uniform bimorphism. The assumption on $\gamma $
to be uniform cannot be removed here. For example, the embedding $\iota :%
\mathbb{N}_{0}\longrightarrow \mathbb{Z}$ is a bimorphism of commutative
monoids ($\mathbb{N}_{0}$-semimodules) which is not an isomorphism. Notice
that $\iota $ is not $i$-uniform; in fact $\overline{\iota (\mathbb{N}_{0}})=%
\mathbb{Z}.$
\end{remark}

\begin{lemma}
\label{1st-IT}\emph{(Compare with \cite[Proposition 4.3.]{Tak1981})} Let $%
\gamma :X\rightarrow Y$ be a morphism of $S$-semimodules.

\begin{enumerate}
\item The sequence%
\begin{equation}
0\rightarrow \mathrm{Ker}(\gamma )\overset{\mathrm{\ker }(\gamma )}{%
\longrightarrow }X\overset{\gamma }{\rightarrow }Y\overset{\mathrm{\mathrm{%
coker}}(\gamma )}{\longrightarrow }\mathrm{\mathrm{\mathrm{\mathrm{\mathrm{%
Coker}}}}}(\gamma )\rightarrow 0  \label{ker-coker}
\end{equation}

is semi-exact. Moreover, (\ref{ker-coker}) is exact if and only if $\gamma $
is uniform.

\item We have two exact sequences%
\begin{equation*}
0\rightarrow \overline{\gamma (X)}\overset{\mathrm{ker}(\mathrm{\mathrm{coker%
}}(\gamma ))}{\longrightarrow }Y\overset{\mathrm{\mathrm{coker}}(\gamma )}{%
\longrightarrow }Y/\gamma (X)\rightarrow 0.
\end{equation*}%
and%
\begin{equation*}
0\rightarrow \mathrm{Ker}(\gamma )\overset{\mathrm{ker}(\gamma )}{%
\longrightarrow }X\overset{\mathrm{\mathrm{coker}}(\mathrm{ker}(\gamma ))}{%
\longrightarrow }X/\mathrm{Ker}(\gamma )\rightarrow 0.
\end{equation*}
\end{enumerate}
\end{lemma}

\begin{corollary}
\label{reg-sub}\emph{(Compare with \cite[Proposition 4.8.]{Tak1981}) }Let $M$
be an $S$-semimodule.

\begin{enumerate}
\item Let $\rho $ an $S$-congruence relation on $M$ and consider the
sequence of $S$-semimodules%
\begin{equation*}
0\longrightarrow \mathrm{Ker}(\pi _{\rho })\overset{\iota _{\rho }}{%
\longrightarrow }M\overset{\rho }{\longrightarrow }M/\rho \longrightarrow 0.
\end{equation*}

\begin{enumerate}
\item $0\rightarrow \mathrm{Ker}(\pi _{\rho })\overset{\iota _{\rho }}{%
\longrightarrow }M\overset{\pi _{\rho }}{\longrightarrow }M/\rho \rightarrow
0$ is exact.

\item $M/\rho =\mathrm{Coker}(\iota _{\rho }),$ whence $M/\rho $ is a normal
quotient.
\end{enumerate}

\item Let $L\leq _{S}M$ an $S$-subsemimodule.

\begin{enumerate}
\item The sequence $0\rightarrow L\overset{\iota }{\longrightarrow }M\overset%
{\pi _{L}}{\longrightarrow }M/L\rightarrow 0$ is semi-exact.

\item $0\rightarrow \overline{L}\overset{\iota }{\longrightarrow }M\overset{%
\pi _{L}}{\longrightarrow }M/L\rightarrow 0$ is exact.

\item The following are equivalent:

\begin{enumerate}
\item $0\rightarrow L\overset{\iota }{\longrightarrow }M\overset{\pi _{L}}{%
\longrightarrow }M/L\rightarrow 0$ is exact;

\item $L\simeq \mathrm{Ker}(\pi _{L});$

\item $0\longrightarrow L\overset{\iota }{\longrightarrow }\overline{L}%
\longrightarrow 0$ is exact;

\item $L$ is a uniform subsemimodule;

\item $L$ is a normal subsemimodule.
\end{enumerate}
\end{enumerate}
\end{enumerate}
\end{corollary}

\section{Homological lemmas}

\qquad In this section we prove some elementary diagram lemma for
semimodules over semirings. These apply in particular to commutative
monoids, considered as semimodules over the semiring of non-negative
integers. Recall that a sequence $A\overset{f}{\longrightarrow }B\overset{g}{%
\longrightarrow }C$ of semimodules is exact iff $f(A)=\mathrm{Ker}(g)$ and $%
g $ is is $k$-uniform (equivalently, $f(A)=\mathrm{Ker}(g)$ and $%
g(b)=g(b^{\prime })\Longrightarrow b+f(a)=b^{\prime }+f(a^{\prime })$ for
some $a,a^{\prime }\in A$).

\qquad The following result can be easily proved using \emph{diagram chasing}
(compare \textquotedblleft 2\textquotedblright\ with \cite[Lemma 1.9]%
{Pat2006}).

\begin{lemma}
\label{short}Consider the following commutative diagram of semimodules%
\begin{equation*}
\xymatrix{ & & 0 \ar[d] \\ L_1 \ar[r]^{f_1} \ar[d]_{\alpha_1} & M_1
\ar[r]^{g_1} \ar[d]_{\alpha_2} & N_1 \ar[d]_{\alpha_3} \\ L_2 \ar[r]^{f_2}
\ar[d] & M_2 \ar[r]^{g_2} & N_2 \\ 0 & & }
\end{equation*}%
and assume that the first and the third columns are exact (i.e. $\alpha _{1}$
is surjective and $\alpha _{3}$ is injective).

\begin{enumerate}
\item Let $\alpha _{2}$ be surjective. If the first row is exact, then the
second row is exact.

\item Let $\alpha _{2}$ be injective. If the second row is exact, then the
first row is exact.

\item Let $a_{2}$ be an isomorphism. The first row is exact if and only if
the second row is exact.
\end{enumerate}
\end{lemma}

\begin{Beweis}
\begin{enumerate}
\item Let $\alpha _{2}$ be surjective and assume that the first row is exact.

\begin{itemize}
\item $f_{2}(L_{2})=\mathrm{Ker}(g_{2}).$

Notice that $g_{2}\circ f_{2}\circ \alpha _{1}=g_{2}\circ \alpha _{2}\circ
f_{1}=\alpha _{3}\circ g_{1}\circ f_{1}=0.$ Since $\alpha _{1}$ is an
epimorphism, we conclude that $g_{2}\circ f_{2}=0;$ in particular, $%
f_{2}(L_{2})\subseteq \mathrm{Ker}(g_{2}).$ On the other hand, let $m_{2}\in
\mathrm{Ker}(g_{2}).$ Since $\alpha _{2}$ is surjective, there exists $%
m_{1}\in M_{1}$ such that $\alpha _{2}(m_{1})=m_{2}.$ Since $\alpha _{3}$ is
a semi-monomorphism and $(\alpha _{3}\circ g_{1})(m_{1})=(g_{2}\circ \alpha
_{2})(m_{1})=0,$ we conclude that $g_{1}(m_{1})=0.$ Since the first row is
exact, there exists $l_{1}\in L_{1}$ such that $m_{1}=f_{1}(l_{1}).$ It
follows that $m_{2}=(\alpha _{2}\circ f_{1})(l_{1})=f_{2}(\alpha
_{1}(l_{1}))\in f_{2}(L_{2}).$

\item $g_{2}$ is $k$-uniform.

Suppose that $g_{2}(m_{2})=g_{2}(m_{2}^{\prime }).$ Since $\alpha _{2}$ is
surjective, there exist $m_{1},m_{1}^{\prime }\in M_{1}$ such that $\alpha
_{2}(m_{1})=m_{2}$ and $\alpha _{2}(m_{1}^{\prime })=m_{2}^{\prime }.$ Since
$\alpha _{3}$ is injective and $(\alpha _{3}\circ g_{1})(m_{1})=(g_{2}\circ
\alpha _{2})(m_{1})=(g_{2}\circ \alpha _{2})(m_{1}^{\prime })=(\alpha
_{3}\circ g_{1})(m_{1}^{\prime })$ we have $g_{1}(m_{1})=g_{1}(m_{1}^{\prime
}).$ Since $g_{1}$ is $k$-uniform and $f_{1}(L_{1})=\mathrm{Ker}(g_{1})$
there exist $l_{1},l_{1}^{\prime }\in L_{1}$ such that $%
m_{1}+f_{1}(l_{1})=m_{1}^{\prime }+f_{1}(l_{1}^{\prime }).$ It follows that $%
m_{2}+(\alpha _{2}\circ f_{1})(l_{1})=m_{2}^{\prime }+(\alpha _{2}\circ
f_{1})(l_{1}^{\prime })$ whence $m_{2}+f_{2}(\alpha
_{1}(l_{1}))=m_{2}^{\prime }+f_{2}(\alpha _{1}(l_{1}^{\prime })).$ Since $%
f_{2}(L_{2})\subseteq \mathrm{Ker}(g_{2}),$ we conclude that $g_{2}$ is $k$%
-uniform.
\end{itemize}

\item Let $\alpha _{2}$ be injective and assume that the second row is exact.

\begin{itemize}
\item $f_{1}(L_{1})=\mathrm{Ker}(g_{1}).$

Notice that $\alpha _{3}\circ g_{1}\circ f_{1}=g_{2}\circ \alpha _{2}\circ
f_{1}=g_{2}\circ f_{2}\circ \alpha _{1}=0.$ Since $\alpha _{3}$ is a
monomorphism, we conclude that $g_{1}\circ f_{1}=0,$ i.e. $%
f_{1}(L_{1})\subseteq \mathrm{Ker}(g_{1}).$ Let $m_{1}\in \mathrm{Ker}%
(g_{1}).$ Then $(g_{2}\circ \alpha _{2})(m_{1})=(\alpha _{3}\circ
g_{1})(m_{1})=0.$ Since the second row is exact, there exist $l_{2}\in L_{2}$
such that $f_{2}(l_{2})=\alpha _{2}(m_{1}).$ Since $\alpha _{1}$ is
surjective, there exists $l_{1}\in L_{1}$ such that $\alpha
_{2}(m_{1})=f_{2}(l_{2})=f_{2}(\alpha _{1}(l_{1}))=(\alpha _{2}\circ
f_{1})(l_{1}).$ Since $\alpha _{2}$ is injective, $m_{1}=f_{1}(l_{1}).$

\item $g_{1}$ is $k$-uniform.

Suppose that $g_{1}(m_{1})=g_{1}(m_{1}^{\prime })$ for some $%
m_{1},m_{1}^{\prime }\in M_{1}.$ Then we have $(g_{2}\circ \alpha
_{2})(m_{1})=(\alpha _{3}\circ g_{1})(m_{1})=(\alpha _{3}\circ
g_{1})(m_{1}^{\prime })=(g_{2}\circ \alpha _{2})(m_{1}^{\prime }).$ Since $%
g_{2}$ is $k$-uniform and $f_{2}(L_{2})=\mathrm{Ker}(g_{2}),$ there exist $%
l_{2},l_{2}^{\prime }\in L_{2}$ such that $\alpha
_{2}(m_{1})+f_{2}(l_{2})=\alpha _{2}(m_{1}^{\prime })+f_{2}(l_{2}^{\prime
}). $ Since $\alpha _{1}$ is surjective, there exist $l_{1},l_{1}^{\prime
}\in L_{1}$ such that $\alpha _{2}(m_{1}+f_{1}(l_{1}))=\alpha
_{2}(m_{1})+(f_{2}\circ \alpha _{1})(l_{1})=\alpha _{2}(m_{1}^{\prime
})+(f_{2}\circ \alpha _{1})(l_{1}^{\prime })=\alpha _{2}(m_{1}^{\prime
}+f_{1}(l_{1}^{\prime })).$ Since $\alpha _{2}$ is injective, we have $%
m_{1}+f_{1}(l_{1})=m_{2}+f_{1}(l_{1}^{\prime })$ and we are done since $%
f_{1}(L_{1})\subseteq \mathrm{Ker}(g_{1}).$
\end{itemize}

\item This is a combination of \textquotedblleft 1\textquotedblright\ and
\textquotedblleft 2\textquotedblright .$\blacksquare $
\end{enumerate}
\end{Beweis}

\subsection*{$\mathcal{R}$-Homological Categories}

\qquad It is well-known that the category of groups, despite being
non-Abelian (in fact not even Puppe-exact, but semiabelian in the sense of
Janelidze et al. \cite{JMT2002}), satisfies the so-called \emph{Short Five
Lemma}. It was shown in \cite[Theorem 4.1.10]{BB2004} that satisfying this
lemma characterizes the so-called \emph{protomodular categories}, whence the
\emph{homological categories}, among the pointed regular ones. Inspired by
this, we introduce in what follows a notion of (\emph{weak}) relative
homological categories with prototype the category of cancellative
commutative monoids, or more generally, the categories of cancellative
semimodules over semirings.

\begin{definition}
Let $\mathfrak{C}$ be a pointed category and $\mathcal{R}=((\mathbf{E},%
\mathbf{M});\mathcal{A})$ where $(\mathbf{E},\mathbf{M})$ is a factorization
structure for $\mathfrak{C}$ and $\mathcal{A}\subseteq \mathrm{Mor}(%
\mathfrak{C}).$ We say that $\mathfrak{C}$ satisfies the \emph{Short }$%
\mathcal{R}$\emph{-Five Lemma} iff for every commutative diagram with $(%
\mathbf{E},\mathbf{M})$-exact rows and $\alpha _{2}\in \mathcal{A}:$%
\begin{equation*}
\xymatrix{0 \ar[r] & L_1 \ar[r]^{f_1} \ar[d]_{\alpha_1} & M_1 \ar[r]^{g_1}
\ar[d]_{\alpha_2} & N_1 \ar[d]_{\alpha_3} \ar[r] & 0\\ 0 \ar[r] & L_2
\ar[r]^{f_2} & M_2 \ar[r]^{g_2} & N_2 \ar[r] & 0}
\end{equation*}%
if $\alpha _{1}$ and $\alpha _{3}$ are isomorphisms, then $\alpha _{2}$ is
an isomorphism.
\end{definition}

\begin{definition}
Let $\mathfrak{C}$ be a category and $\mathcal{R}=((\mathbf{E},\mathbf{M});%
\mathcal{A})$ where $\mathbf{E},\mathbf{M},\mathcal{A}\subseteq \mathrm{Mor}(%
\mathfrak{C}).$ We say that $\mathfrak{C}$ is

\begin{enumerate}
\item $(\mathbf{E},\mathbf{M})$\emph{-regular} iff $\mathfrak{C}$ has finite
limits, is $(\mathbf{E},\mathbf{M})$-structured and the morphisms in $%
\mathbf{E}$ are \emph{pullback stable}.

\item $\mathcal{R}$\emph{-homological category} iff $\mathfrak{C}$ is $(%
\mathbf{E},\mathbf{M})$-regular and satisfies the Short $\mathcal{R}$-Five
Lemma.
\end{enumerate}
\end{definition}

\begin{ex}
One recovers the homological categories in the sense of \cite{BB2004} (i.e.
those which are pointed, regular and protomodular) as follows: a pointed
category $\mathfrak{C}$ is homological iff $\mathfrak{C}$ is $\mathcal{R}$%
-homological where $\mathcal{R}=((\mathbf{RegEpi},\mathbf{Mono});\mathrm{Mor}%
(\mathfrak{C})).$
\end{ex}

\begin{lemma}
\label{diagram}Consider the following commutative diagram of semimodules
with exact rows%
\begin{equation*}
\xymatrix{L_1 \ar[r]^{f_1} \ar[d]_{\alpha_1} & M_1 \ar[r]^{g_1}
\ar[d]_{\alpha_2} & N_1 \ar[d]_{\alpha_3} \\ L_2 \ar[r]^{f_2} & M_2
\ar[r]^{g_2} & N_2}
\end{equation*}

\begin{enumerate}
\item We have:

\begin{enumerate}
\item Let $g_{1}$ and $\alpha _{1}$ be surjective. If $\alpha _{2}$ is
injective, then $\alpha _{3}$ is injective.

\item Let $f_{2}$ be injective and $\alpha _{3}$ a semi-monomorphism. If $%
\alpha _{2}$ is surjective, then $\alpha _{1}$ is surjective.
\end{enumerate}

\item Let $f_{2}$ be a semi-monomorphism.

\begin{enumerate}
\item If $\alpha _{1}$ and $\alpha _{3}$ are semi-monomorphisms, then $%
\alpha _{2}$ is a semi-monomorphism.

\item Let $f_{1},$ $\alpha _{2}$ be cancellative and $f_{2}$ be $k$-uniform.
If $\alpha _{1}$ and $\alpha _{3}$ are injective, then $\alpha _{2}$ is
injective.

\item If $g_{1},$ $\alpha _{1},$ $\alpha _{3}$ are surjective (and $\alpha
_{2}$ is $i$-uniform), then $\alpha _{2}$ is a semi-epimorphism (surjective).
\end{enumerate}
\end{enumerate}
\end{lemma}

\begin{Beweis}
\begin{enumerate}
\item Consider the given commutative diagram.

\begin{enumerate}
\item $\alpha _{3}$ is injective.

Suppose that $\alpha _{3}(n_{1})=\alpha _{3}(n_{1}^{\prime })$ for some $%
n_{1},n_{1}^{\prime }\in N_{1}.$ Since $g_{1}$ is surjective, $%
n_{1}=g_{1}(m_{1})$ and $n_{1}^{\prime }=g_{1}(m_{1}^{\prime })$ for some $%
m_{1},m_{1}^{\prime }\in M_{1}.$ It follows that $(g_{2}\circ \alpha
_{2})(m_{1})=(g_{2}\circ \alpha _{2})(m_{1}^{\prime }).$ Since $g_{2}$ is $k$%
-uniform and $f_{2}(L_{2})=\mathrm{Ker}(g_{2}),$ there exist $%
l_{2},l_{2}^{\prime }\in L_{2}$ such that $\alpha
_{2}(m_{1})+f_{2}(l_{2})=\alpha _{2}(m_{1}^{\prime })+f_{2}(l_{2}^{\prime
}). $ By assumption, $\alpha _{1}$ is surjective and so there exist $%
l_{1},l_{1}^{\prime }\in L_{1}$ such that $\alpha _{1}(l_{1})=l_{2}$ and $%
\alpha _{1}(l_{1}^{\prime })=l_{2}^{\prime }.$ It follows that%
\begin{equation*}
\begin{array}{rclc}
\alpha _{2}(m_{1})+(f_{2}\circ \alpha _{1})(l_{1}) & = & \alpha
_{2}(m_{1}^{\prime })+(f_{2}\circ \alpha _{1})(l_{1}^{\prime }) &  \\
\alpha _{2}(m_{1})+(\alpha _{2}\circ f_{1})(l_{1}) & = & \alpha
_{2}(m_{1}^{\prime })+(\alpha _{2}\circ f_{1})(l_{1}^{\prime }) &  \\
m_{1}+f_{1}(l_{1}) & = & m_{1}^{\prime }+f_{1}(l_{1}^{\prime }) & \text{(}%
\alpha _{2}\text{ is injective)} \\
g_{1}(m_{1}) & = & g_{1}(m_{1}) & \text{(}g_{1}\circ f_{1}=0\text{)} \\
n_{1} & = & n_{1}^{\prime } &
\end{array}%
\end{equation*}

\item $\alpha _{1}$ is surjective.

Let $l_{2}\in L_{2}.$ Since $\alpha _{2}$ is surjective, there exists $%
m_{1}\in M_{1}$ such that $f_{2}(l_{2})=\alpha _{2}(m_{1}).$ It follows that
$0=(g_{2}\circ f_{2})(l_{2})=(g_{2}\circ \alpha _{2})(m_{1})=(\alpha
_{3}\circ g_{1})(m_{1}),$ whence $g_{1}(m_{1})=0$ (since $\alpha _{3}$ is a
semi-monomorphism). Since the first row is exact, $m_{1}=f_{1}(l_{1})$ for
some $l_{1}\in L_{1}$ and so $f_{2}(l_{2})=\alpha _{2}(m_{1})=(\alpha
_{2}\circ f_{1})(l_{1})=(f_{2}\circ \alpha _{1})(l_{1}).$ Since $f_{2}$ is
injective, we have $l_{2}=\alpha _{1}(l_{1}).$
\end{enumerate}

\item Let $f_{2}$ be a semi-monomorphism, \emph{i.e.} $\mathrm{Ker}%
(f_{2})=0. $

\begin{enumerate}
\item We claim that $\alpha _{2}$ is a semi-monomorphism.

Suppose that $\alpha _{2}(m_{1})=0$ for some $m_{1}\in M_{1}.$ Then $(\alpha
_{3}\circ g_{1})(m_{1})=(g_{2}\circ \alpha _{2})(m_{1})=0,$ whence $%
g_{1}(m_{1})=0$ since $\mathrm{Ker}(\alpha _{3})=0.$ It follows that $%
m_{1}=f_{1}(l_{1})$ for some $l_{1}\in L_{1}.$ So, $0=\alpha
_{2}(m_{1})=(\alpha _{2}\circ f_{1})(l_{1})=(f_{2}\circ \alpha _{1})(l_{1}),$
whence $l_{1}=0$ since both $f_{2}$ and $\alpha _{1}$ are
semi-monormorphisms; consequently, $m_{1}=f_{1}(l_{1})=0.$

\item We claim that $\alpha _{2}$ is injective.

Suppose that $\alpha _{2}(m_{1})=\alpha _{2}(m_{1}^{\prime })$ for some $%
m_{1},m_{1}^{\prime }\in M_{1}.$ Then $(\alpha _{3}\circ
g_{1})(m_{1})=(g_{2}\circ \alpha _{2})(m_{1})=(g_{2}\circ \alpha
_{2})(m_{1}^{\prime })=(\alpha _{3}\circ g_{1})(m_{1}^{\prime }),$ whence $%
g_{1}(m_{1})=g_{1}(m_{1}^{\prime })$ since $\alpha _{3}$ is injective. Since
$g_{1}$ is $k$-uniform and $\mathrm{Ker}(g_{1})=f_{1}(L_{1}),$ there exist $%
l_{1},l_{1}^{\prime }\in L_{1}$ such that $m_{1}+f_{1}(l_{1})=m_{1}^{\prime
}+f_{1}(l_{1}^{\prime }).$ Then we have%
\begin{equation*}
\begin{array}{rclc}
\alpha _{2}(m_{1})+(\alpha _{2}\circ f_{1})(l_{1}) & = & \alpha
_{2}(m_{1}^{\prime })+(\alpha _{2}\circ f_{1})(l_{1}^{\prime }) &  \\
\alpha _{2}(m_{1}^{\prime })+(f_{2}\circ \alpha _{1})(l_{1}) & = & \alpha
_{2}(m_{1}^{\prime })+(f_{2}\circ \alpha _{1})(l_{1}^{\prime }) &  \\
(f_{2}\circ \alpha _{1})(l_{1}) & = & (f_{2}\circ \alpha _{1})(l_{1}^{\prime
}) & \text{(}\alpha _{2}\text{ is cancellative)} \\
l_{1} & = & l_{1}^{\prime } & \text{(}f_{2}\text{ and }\alpha _{1}\text{ are
injective)} \\
m_{1}+f_{1}(l_{1}^{\prime }) & = & m_{1}^{\prime }+f_{1}(l_{1}^{\prime }) &
\text{(}f_{1}\text{ is cancellative)} \\
m_{1} & = & m_{1}^{\prime } &
\end{array}%
\end{equation*}

\item We claim that $\alpha _{2}$ is a semi-epimorphism.

Let $m_{2}\in M_{2}.$ Since $\alpha _{3}$ and $g_{1}$ are surjective, there
exists $m_{1}\in M_{1}$ such that $g_{2}(m_{2})=(\alpha _{3}\circ
g_{1})(m_{1})=(g_{2}\circ \alpha _{2})(m_{1}).$ Since $g_{2}$ is $k$%
-uniform, $f_{2}(L_{2})=\mathrm{Ker}(g_{2})$ and $\alpha _{1}$ is
surjective, there exist $l_{1},l_{1}^{\prime }\in L_{1}$ such that%
\begin{eqnarray*}
m_{2}+(f_{2}\circ \alpha _{1})(l_{1}) &=&\alpha _{2}(m_{1})+(f_{2}\circ
\alpha _{1})(l_{1}^{\prime }) \\
m_{2}+\alpha _{2}(f_{1}(l_{1})) &=&\alpha _{2}(m_{1}+f_{1}(l_{1}^{\prime })).
\end{eqnarray*}%
Consequently, $M_{2}=\overline{\alpha _{2}(M_{1})},$ \emph{i.e.} $\alpha
_{2} $ is a semi-epimorphism. If $\alpha _{2}$ is $i$-uniform, then $M_{2}=%
\overline{\alpha _{2}(M_{1})}=\alpha _{2}(M_{1}),$ whence $\alpha _{2}$ is
surjective.$\blacksquare $
\end{enumerate}
\end{enumerate}
\end{Beweis}

\begin{corollary}
\label{cor-short5}Consider the following commutative diagram of semimodules
with exact rows and assume that $M_{1}$ and $M_{2}$ are cancellative%
\begin{equation*}
\xymatrix{& L_1 \ar[r]^{f_1} \ar[d]_{\alpha_1} & M_1 \ar[r]^{g_1}
\ar[d]_{\alpha_2} & N_1 \ar[d]_{\alpha_3} \ar[r] & 0\\ 0 \ar[r] & L_2
\ar[r]^{f_2} & M_2 \ar[r]^{g_2} & N_2 & }
\end{equation*}

\begin{enumerate}
\item Let $\alpha _{2}$ be an isomorphism. Then $\alpha _{1}$ is surjective
if and only if $\alpha _{3}$ is injective.

\item Let $\alpha _{2}$ be $i$-uniform. If $\alpha _{1}$ and $\alpha _{3}$
are isomorphisms, then $\alpha _{2}$ is an isomorphism.
\end{enumerate}
\end{corollary}

\begin{proposition}
\label{short-5}\emph{(The Short Five Lemma)} Consider the following
commutative diagram of semimodules with $M_{1},M_{2}$ cancellative%
\begin{equation*}
\xymatrix{0 \ar[r] & L_1 \ar[r]^{f_1} \ar[d]_{\alpha_1} & M_1 \ar[r]^{g_1}
\ar[d]_{\alpha_2} & N_1 \ar[d]_{\alpha_3} \ar[r] & 0\\ 0 \ar[r] & L_2
\ar[r]^{f_2} & M_2 \ar[r]^{g_2} & N_2 \ar[r] & 0}
\end{equation*}%
Then $\alpha _{1},$ $\alpha _{3}$ are isomorphisms and $\alpha _{2}$ is $i$%
-uniform if and only if $\alpha _{2}$ is an isomorphism. In particular, the
category $\mathbb{CS}_{S}$ of cancellative right $S$-semimodules is $%
\mathcal{R}$-homological, where $\mathcal{R}=((\mathbf{Surj},\mathbf{Inj});%
\mathcal{I})$ and $\mathcal{I}$ is the class of $i$-uniform morphisms.
\end{proposition}

\begin{lemma}
\label{5-details}Consider the following commutative diagram of semimodules
with exact rows%
\begin{equation*}
\xymatrix{U_1 \ar[r]^{e_1} \ar[d]_{\gamma} & L_1 \ar[r]^{f_1}
\ar[d]_{\alpha_1} & M_1 \ar[r]^{g_1} \ar[d]_{\alpha_2} & N_1 \ar[r]
\ar[d]_{\alpha_3} \ar[r]^{h_1} & V_1 \ar[d]_{\delta} \\ U_2 \ar[r]^{e_2} &
L_2 \ar[r]^{f_2} & M_2 \ar[r]^{g_2} & N_2 \ar[r]^{h_2} & V_2}
\end{equation*}

\begin{enumerate}
\item Let $\gamma $ be surjective.

\begin{enumerate}
\item If $\alpha _{1}$ is injective and $\alpha _{3}$ is a
semi-monomorphisms, then $\alpha _{2}$ is a semi-monomorphism.

\item Assume that $f_{1}$ and $\alpha _{2}$ are cancellative. If $\alpha
_{1} $ and $\alpha _{3}$ are injective, then $\alpha _{2}$ is injective.
\end{enumerate}

\item Let $\delta $ be a semi-monomorphism. If $\alpha _{1},$ $\alpha _{3}$
are surjective (and $\alpha _{2}$ is $i$-uniform), then $\alpha _{2}$ is a
semi-epimorphism (surjective).

\item Let $f_{1},\alpha _{2}$ be cancellative, $\gamma $ be surjective and $%
\delta $ be injective. If $\alpha _{1}$ and $\alpha _{3}$ are isomorphisms,
then $\alpha _{2}$ is injective and a semi-epimorphism.
\end{enumerate}
\end{lemma}

\begin{Beweis}
Assume that the diagram is commutative and that the two rows are exact.

\begin{enumerate}
\item Let $\gamma $ be surjective.

\begin{enumerate}
\item Assume that $\alpha _{1}$ is injective and that $\alpha _{3}$ is a
semi-isomorphism. We claim that $\alpha _{2}$ is a semi-monomorphism.

Suppose that $\alpha _{2}(m_{1})=0$ for some $m_{1}\in M_{1}$ so that $%
(\alpha _{3}\circ g_{1})(m_{1})=(g_{2}\circ \alpha _{2})(m_{1})=0.$ Since $%
\alpha _{3}$ is a semi-monomorphism $g_{1}(m_{1})=0,$ whence $%
m_{1}=f_{1}(l_{1})$ for some $l_{1}\in L_{1}.$ So, $0=\alpha
_{2}(m_{1})=(\alpha _{2}\circ f_{1})(l_{1})=(f_{2}\circ \alpha _{1})(l_{1}),$
whence $\alpha _{1}(l_{1})=(e_{2}\circ \gamma )(u_{1})=(\alpha _{1}\circ
e_{1})(u_{1})$ for some $u_{1}\in U_{1}$ (since $\gamma $ is surjective and $%
\mathrm{Ker}(f_{2})=e_{2}(U_{2})$). Since $\alpha _{1}$ is injective, it
follows that $l_{1}=e_{1}(u_{1})$ whence $m_{1}=f_{1}(l_{1})=(f_{1}\circ
e_{1})(u_{1})=0.$

\item Assume that $f_{1},\alpha _{2}$ are cancellative and $\alpha _{1},$ $%
\alpha _{3}$ are injective. We claim that $\alpha _{2}$ is injective.

Suppose that $\alpha _{2}(m_{1})=\alpha _{2}(m_{1}^{\prime })$ for some $%
m_{1},m_{1}^{\prime }\in M_{1}.$ Then $(\alpha _{3}\circ
g_{1})(m_{1})=(g_{2}\circ \alpha _{2})(m_{1})=(g_{2}\circ \alpha
_{2})(m_{1}^{\prime })=(\alpha _{3}\circ g_{1})(m_{1}^{\prime }),$ whence $%
g_{1}(m_{1})=g_{1}(m_{1}^{\prime })$ (notice that $\alpha _{3}$ is
injective). Since $g_{1}$ is $k$-uniform and $\mathrm{Ker}%
(g_{1})=f_{1}(L_{1}),$ there exist $l_{1},l_{1}^{\prime }\in L_{1}$ such
that $m_{1}+f_{1}(l_{1})=m_{1}^{\prime }+f_{1}(l_{1}^{\prime }).$ Then we
have%
\begin{equation*}
\begin{array}{rclc}
\alpha _{2}(m_{1})+(\alpha _{2}\circ f_{1})(l_{1}) & = & \alpha
_{2}(m_{1}^{\prime })+(\alpha _{2}\circ f_{1})(l_{1}^{\prime }) &  \\
\alpha _{2}(m_{1}^{\prime })+(f_{2}\circ \alpha _{1})(l_{1}) & = & \alpha
_{2}(m_{1}^{\prime })+(f_{2}\circ \alpha _{1})(l_{1}^{\prime }) &  \\
f_{2}(\alpha _{1}(l_{1})) & = & f_{2}(\alpha _{1}(l_{1}^{\prime })) & \text{(%
}\alpha _{2}\text{ is cancellative)} \\
\alpha _{1}(l_{1})+k_{2} & = & \alpha _{1}(l_{1}^{\prime })+k_{2}^{\prime }
& \text{(}f_{2}\text{ is $k$-uniform)} \\
\alpha _{1}(l_{1})+(e_{2}\circ \gamma )(u_{1}) & = & \alpha
_{1}(l_{1}^{\prime })+(e_{2}\circ \gamma )(u_{1}^{\prime }) & \text{(}\gamma
\text{ is surjective)} \\
\alpha _{1}(l_{1})+(\alpha _{1}\circ e_{1})(u_{1}) & = & \alpha
_{1}(l_{1}^{\prime })+(\alpha _{1}\circ e_{1})(u_{1}^{\prime }) &  \\
l_{1}+e_{1}(u_{1}) & = & l_{1}^{\prime }+e_{1}(u_{1}^{\prime }) & \text{(}%
\alpha _{1}\text{ is injective)} \\
f_{1}(l_{1}) & = & f_{1}(l_{1}^{\prime }) & \text{(}f_{1}\circ e_{1}=0\text{)%
} \\
m_{1}+f_{1}(l_{1}) & = & m_{1}+f_{1}(l_{1}^{\prime }) &  \\
m_{1}^{\prime }+f_{1}(l_{1}^{\prime }) & = & m_{1}+f_{1}(l_{1}^{\prime }) &
\\
m_{1}^{\prime } & = & m_{1} & \text{(}f_{1}\text{ is cancellative)}%
\end{array}%
\end{equation*}
\end{enumerate}

\item Let $\delta $ be a semi-monomorphism. Assume that $\alpha _{1}$ and $%
\alpha _{3}$ are surjective. Let $m_{2}\in M_{2}.$ Since $\alpha _{3}$ is
surjective, there exists $n_{1}\in N_{1}$ such that $g_{2}(m_{2})=\alpha
_{3}(n_{1}).$ It follows that $0=(h_{2}\circ g_{2})(m_{2})=(h_{2}\circ
\alpha _{3})(n_{1})=(\delta \circ h_{1})(n_{1}),$ whence $h_{1}(n_{1})=0$
since $\delta $ is a semi-monomorphism. Since $g_{1}(M_{1})=\mathrm{Ker}%
(h_{1}),$ we have $n_{1}=g_{1}(m_{1})$ for some $m_{1}\in M_{1}.$ Notice
that $(g_{2}\circ \alpha _{2})(m_{1})=(\alpha _{3}\circ g_{1})(m_{1})=\alpha
_{3}(n_{1})=g_{2}(m_{2}).$ Since $g_{2}$ is $k$-uniform, $f_{2}(L_{2})=%
\mathrm{Ker}(g_{2})$ and $\alpha _{1}$ is surjective, there exists $%
l_{1},l_{1}^{\prime }\in L_{1}$ such that%
\begin{eqnarray*}
m_{2}+(f_{2}\circ \alpha _{1})(l_{1}) &=&\alpha _{2}(m_{1})+(f_{2}\circ
\alpha _{1})(l_{1}^{\prime }) \\
m_{2}+\alpha _{2}(f_{1}(l_{1})) &=&\alpha _{2}(m_{1}+f_{1}(l_{1}^{\prime })),
\end{eqnarray*}%
i.e. $m_{2}\in \overline{\alpha _{2}(M_{1})}.$ Consequently, $M_{2}=%
\overline{\alpha _{2}(M_{1})}.$ If $\alpha _{2}$ is $i$-uniform, then $%
\alpha _{2}(M)=\overline{\alpha _{2}(M_{1})}=M_{2},$ \emph{i.e. }$\alpha
_{2} $ is surjective.

\item This is a combination of \textquotedblleft 1\textquotedblright\ and
\textquotedblleft 2\textquotedblright .$\blacksquare $
\end{enumerate}
\end{Beweis}

\begin{corollary}
\label{5-lemma}\emph{(The Five Lemma) }Consider the following commutative
diagram of semimodules with exact rows and columns and assume that $f_{1}$
and $\alpha _{2}$ are cancellative%
\begin{equation*}
\xymatrix{ & & & & 0 \ar[d] \\ U_1 \ar[r]^{e_1} \ar[d]_{\gamma} & L_1
\ar[r]^{f_1} \ar[d]_{\alpha_1} & M_1 \ar[r]^{g_1} \ar[d]_{\alpha_2} & N_1
\ar[r] \ar[d]_{\alpha_3} \ar[r]^{h_1} & V_1 \ar[d]_{\delta} \\ U_2
\ar[r]^{e_2} \ar[d] & L_2 \ar[r]^{f_2} & M_2 \ar[r]^{g_2} & N_2 \ar[r]^{h_2}
& V_2 \\ 0 & & & & }
\end{equation*}

\begin{enumerate}
\item If $\alpha _{1}$ and $\alpha _{3}$ are injective, then $\alpha _{2}$
is injective.

\item Let $\alpha _{2}$ be $i$-uniform. If $\alpha _{1}$ and $\alpha _{3}$
are surjective, then $\alpha _{2}$ is surjective.

\item Let $\alpha _{2}$ be $i$-uniform. If $\alpha _{1}$ and $\alpha _{3}$
are isomorphisms, then $\alpha _{2}$ is an isomorphism.
\end{enumerate}
\end{corollary}

\subsection*{The Snake Lemma}

\qquad One of the basic homological lemmas that are proved usually in
categories of modules (e.g. \cite{Wis1991}), or more generally in Abelian
categories, is the so called \emph{Kernel-Cokernel Lemma} (\emph{Snake Lemma}%
). Several versions of this lemma were proved also in non-abelian categories
(e.g. \emph{homological categories} \cite{BB2004}, \emph{relative
homological categories} \cite{Jan2006} and incomplete relative homological
categories \cite{Jan2010b}).

\begin{lemma}
\label{9-1}Consider the following commutative diagram with exact columns and
assume that the second row is exact.%
\begin{equation*}
\xymatrix{ & & 0 \ar[d] & 0 \ar[d] & \\ & L_1 \ar[r]^{f_1} \ar[d]_{\alpha_1}
& M_1 \ar[r]^{g_1} \ar[d]_{\alpha_2} & N_1 \ar[d]_{\alpha_3} & \\ & L_2
\ar[r]^{f_2} \ar[d]_{\beta_1} & M_2 \ar[r]^{g_2} \ar[d]_{\beta_2} & N_2
\ar[d]_{\beta_3} &\\ & L_3 \ar[r]^{f_3} & M_3 \ar[r]^{g_3} & N_3 & }
\end{equation*}

\begin{enumerate}
\item If $f_{3}$ is injective and $f_{2}$ is cancellative, then the first
row is exact.

\item If $g_{2},$ $\beta _{1}$ are surjective, the third row is exact (and $%
g_{1}$ is $i$-uniform), then $g_{1}$ is a semi-epimorphism (surjective).
\end{enumerate}
\end{lemma}

\begin{Beweis}
Assume that the second row is exact.

\begin{enumerate}
\item Notice that $\alpha _{3}\circ g_{1}\circ f_{1}=g_{2}\circ \alpha
_{2}\circ f_{1}=g_{2}\circ f_{2}\circ \alpha _{1}=0,$ whence $g_{1}\circ
f_{1}=0$ since $\alpha _{3}$ is a monomorphism. In particular, $%
f_{1}(L_{1})\subseteq \mathrm{Ker}(g_{1}).$

\begin{itemize}
\item We claim that $f_{1}(L_{1})=\mathrm{Ker}(g_{1}).$

Let $m_{1}\in \mathrm{Ker}(g_{1}),$ so that $g_{1}(m_{1})=0.$ It follows that%
\begin{equation*}
\begin{array}{rclc}
(\alpha _{3}\circ g_{1})(m_{1}) & = & 0 &  \\
(g_{2}\circ \alpha _{2})(m_{1}) & = & 0 &  \\
\alpha _{2}(m_{1}) & = & f_{2}(l_{2}) & \text{(2nd row is proper exact)} \\
0 & = & (\beta _{2}\circ f_{2})(l_{2}) & \text{(}\beta _{2}\circ \alpha
_{2}=0\text{)} \\
0 & = & (f_{3}\circ \beta _{1})(l_{2}) &  \\
\beta _{1}(l_{2}) & = & 0 & \text{(}f_{3}\text{ is a semi-monomorphism)} \\
l_{2} & = & \alpha _{1}(l_{1}) & \text{(1st column is proper exact)} \\
f_{2}(l_{2}) & = & (f_{2}\circ \alpha _{1})(l_{1}) &  \\
\alpha _{2}(m_{1}) & = & \alpha _{2}(f_{1}(l_{1})) &  \\
m_{1} & = & f_{1}(l_{1}) & \text{(}\alpha _{2}\text{ is injective)}%
\end{array}%
\end{equation*}

\item We claim that $g_{1}$ is $k$-uniform.

Suppose that $g_{1}(m_{1})=g_{1}(m_{1}^{\prime })$ for some $%
m_{1},m_{1}^{\prime }\in M_{1}.$ It follows that%
\begin{equation*}
\begin{array}{rclc}
(\alpha _{3}\circ g_{1})(m_{1}) & = & (\alpha _{3}\circ g_{1})(m_{1}^{\prime
}) &  \\
(g_{2}\circ \alpha _{2})(m_{1}) & = & (g_{2}\circ \alpha _{2})(m_{1}^{\prime
}) &  \\
\alpha _{2}(m_{1})+f_{2}(l_{2}) & = & \alpha _{2}(m_{1}^{\prime
})+f_{2}(l_{2}^{\prime })\text{ (2nd row is exact)} &  \\
(\beta _{2}\circ f_{2})(l_{2}) & = & (\beta _{2}\circ f_{2})(l_{2}^{\prime })%
\text{ (}\beta _{2}\circ \alpha _{2}=0\text{)} &  \\
(f_{3}\circ \beta _{1})(l_{2}) & = & (f_{3}\circ \beta _{1})(l_{2}^{\prime })
&  \\
\beta _{1}(l_{2}) & = & \beta _{1}(l_{2}^{\prime })\text{ (}f_{3}\text{ is
injective)} &  \\
l_{2}+\alpha _{1}(l_{1}) & = & l_{2}^{\prime }+\alpha _{1}(l_{1}^{\prime })%
\text{ (first column is exact)} &  \\
f_{2}(l_{2})+(f_{2}\circ \alpha _{1})(l_{1}) & = & f_{2}(l_{2}^{\prime
})+(f_{2}\circ \alpha _{1})(l_{1}^{\prime }) &  \\
f_{2}(l_{2})+(\alpha _{2}\circ f_{1})(l_{1}) & = & f_{2}(l_{2}^{\prime
})+(\alpha _{2}\circ f_{1})(l_{1}^{\prime }) &  \\
\alpha _{2}(m_{1})+f_{2}(l_{2})+(\alpha _{2}\circ f_{1})(l_{1}) & = & \alpha
_{2}(m_{1})+f_{2}(l_{2}^{\prime })+(\alpha _{2}\circ f_{1})(l_{1}^{\prime })
&  \\
f_{2}(l_{2}^{\prime })+\alpha _{2}(m_{1}^{\prime }+f_{1}(l_{1})) & = &
f_{2}(l_{2}^{\prime })+\alpha _{2}(m_{1}+f_{1}(l_{1}^{\prime }))\text{ (}%
f_{2}\text{ is cancellative)} &  \\
m_{1}^{\prime }+f_{1}(l_{1}) & = & m_{1}+f_{1}(l_{1}^{\prime })\text{ (}%
\alpha _{2}\text{ is injective)} &
\end{array}%
\end{equation*}%
Since $f_{1}(L_{1})\subseteq \mathrm{Ker}(g_{1}),$ it follows that $g_{1}$
is $k$-uniform.
\end{itemize}

\item We claim that $g_{1}$ is a semi-epimorphism.

Let $n_{1}\in N_{1}.$ Let $m_{2}\in M_{2}$ be such that $g_{2}(m_{2})=\alpha
_{3}(n_{1}).$ Then%
\begin{equation*}
\begin{array}{rclc}
g_{3}(\beta _{2}(m_{2})) & = & \beta _{3}(g_{2}(m_{2})) &  \\
& = & (\beta _{3}\circ \alpha _{3})(m_{2}) &  \\
& = & 0 & \text{(}\beta _{3}\circ \alpha _{3}=0\text{)} \\
\beta _{2}(m_{2}) & = & f_{3}(l_{3}) & \text{(3rd row is exact)} \\
& = & f_{3}(\beta _{1}(l_{2})) & \text{(}\beta _{1}\text{ is surjective)} \\
& = & \beta _{2}(f_{2}(l_{2})) &  \\
m_{2}+\alpha _{2}(m_{1}) & = & f_{2}(l_{2})+\alpha _{2}(m_{1}^{\prime }) &
\text{(2nd column is exact)} \\
g_{2}(m_{2})+(g_{2}\circ \alpha _{2})(m_{1}) & = & (g_{2}\circ \alpha
_{2})(m_{1}) & \text{(}g_{2}\circ f_{2}=0\text{)} \\
\alpha _{3}(n_{1}+g_{1}(m_{1})) & = & \alpha _{3}(g_{1}(m_{1}^{\prime })) &
\\
n_{1}+g_{1}(m_{1}) & = & g_{1}(m_{1}^{\prime }) & \text{(}\alpha _{3}\text{
is injective)}%
\end{array}%
\end{equation*}%
Consequently, $N_{1}=\overline{g_{1}(M_{1})}$ ($=$ $g_{1}(M_{1})$ if $g_{1}$
is assumed to be $i$-uniform).$\blacksquare $
\end{enumerate}
\end{Beweis}

\qquad Similarly, one can prove the following result.

\begin{lemma}
\label{9-3}Consider the following commutative diagram with exact columns and
assume that the second row is exact%
\begin{equation*}
\xymatrix{ & L_1 \ar[r]^{f_1} \ar[d]_{\alpha_1} & M_1 \ar[r]^{g_1}
\ar[d]_{\alpha_2} & N_1 \ar[d]_{\alpha_3} & \\& L_2 \ar[r]^{f_2}
\ar[d]_{\beta_1} & M_2 \ar[r]^{g_2} \ar[d]_{\beta_2} & N_2 \ar[d]_{\beta_3}
&\\ & L_3 \ar[r]^{f_3} \ar[d] & M_3 \ar[r]^{g_3} \ar[d] & N_3 & \\ & 0 & 0 &
}
\end{equation*}

\begin{enumerate}
\item If $g_{1}$ is surjective and $f_{3}$ is $i$-uniform, then the third
row is exact.

\item If $f_{2},$ $\alpha _{3}$ are injective, $\alpha _{2}$ is cancellative
and the first row is exact, then $f_{3}$ is injective.
\end{enumerate}
\end{lemma}

\begin{Beweis}
Assume that the second row is exact.

\begin{enumerate}
\item Notice that $g_{3}\circ f_{3}\circ \beta _{1}=g_{3}\circ \beta
_{2}\circ f_{2}=\beta _{3}\circ g_{2}\circ f_{2}=0.$ Since $\beta _{1}$ is
an epimorphism, we have $g_{3}\circ f_{3}=0$ (i.e. $f_{3}(L_{3})\subseteq
\mathrm{Ker}(g_{3})$).

\begin{itemize}
\item We claim that $f_{3}(L_{3})=\mathrm{Ker}(g_{3}).$ Let $m_{3}\in
\mathrm{Ker}(g_{3}).$

Since $\beta _{2}$ is surjective, $m_{3}=\beta _{2}(m_{2})$ for some $%
m_{2}\in M_{2}.$ It follows that $0=(g_{3}\circ \beta _{2})(m_{2})=(\beta
_{3}\circ g_{2})(m_{2}),$ i.e. $g_{2}(m_{2})\in \mathrm{Ker}(\beta
_{3})=\alpha _{3}(N_{1}).$ We have%
\begin{equation*}
\begin{array}{rclc}
g_{2}(m_{2}) & = & \alpha _{3}(n_{1}) &  \\
& = & (\alpha _{3}\circ g_{1})(m_{1}) & \text{(}g_{1}\text{ is surjective)}
\\
& = & (g_{2}\circ \alpha _{2})(m_{1}) &  \\
m_{2}+f_{2}(l_{2}) & = & \alpha _{2}(m_{1})+f_{2}(l_{2}^{\prime }) & \text{%
(2nd row is exact)} \\
\beta _{2}(m_{2})+(\beta _{2}\circ f_{2})(l_{2}) & = & (\beta _{2}\circ
f_{2})(l_{2}^{\prime }) & \text{(}\beta _{2}\circ \alpha _{2}=0\text{)} \\
m_{3}+(f_{3}\circ \beta _{1})(l_{2}) & = & (f_{3}\circ \beta
_{1})(l_{2}^{\prime }) &
\end{array}%
\end{equation*}%
We conclude that $\mathrm{Ker}(g_{3})=\overline{f_{3}(L_{3})}=f_{3}(L_{3}).$

\item We claim that $g_{3}$ is $k$-uniform.

Suppose that $g_{3}(m_{3})=g_{3}(m_{3}^{\prime })$ for some $%
m_{3},m_{3}^{\prime }\in M_{3}.$ Since $\beta _{2}$ is surjective, there
exist $m_{2},m_{2}^{\prime }\in M$ such that $\beta _{2}(m_{2})=m_{3}$ and $%
\beta _{2}(m_{2}^{\prime })=m_{3}^{\prime }.$ Then%
\begin{equation*}
\begin{array}{rclc}
(g_{3}\circ \beta _{2})(m_{2}) & = & (g_{3}\circ \beta _{2})(m_{2}^{\prime })
&  \\
(\beta _{3}\circ g_{2})(m_{2}) & = & (\beta _{3}\circ g_{2})(m_{2}^{\prime })
&  \\
g_{2}(m_{2})+\alpha _{3}(n_{1}) & = & g_{2}(m_{2}^{\prime })+\alpha
_{3}(n_{1}^{\prime }) & \text{(3rd column is exact)} \\
g_{2}(m_{2})+(\alpha _{3}\circ g_{1})(m_{1}) & = & g_{2}(m_{2}^{\prime
})+(\alpha _{3}\circ g_{1})(m_{1}^{\prime }) & \text{(}g_{1}\text{ is
surjective)} \\
g_{2}(m_{2})+(g_{2}\circ \alpha _{2})(m_{1}) & = & g_{2}(m_{2}^{\prime
})+(g_{2}\circ \alpha _{2})(m_{1}^{\prime }) &  \\
m_{2}+\alpha _{2}(m_{1})+f_{2}(l_{2}) & = & m_{2}^{\prime }+\alpha
_{2}(m_{1}^{\prime })+f_{2}(l_{2}^{\prime }) & \text{(2nd row is exact)} \\
\beta _{2}(m_{2})+(\beta _{2}\circ f_{2})(l_{2}) & = & \beta
_{2}(m_{2}^{\prime })+(\beta _{2}\circ f_{2})(l_{2}^{\prime }) & \text{(}%
\beta _{2}\circ \alpha _{2}=0\text{)} \\
m_{3}+(f_{3}\circ \beta _{1})(l_{2}) & = & m_{3}^{\prime }+(f_{3}\circ \beta
_{1})(l_{2}^{\prime }) &
\end{array}%
\end{equation*}%
Since $f_{3}(L_{3})\subseteq \mathrm{Ker}(g_{3}),$ we conclude that $g_{3}$
is $k$-uniform.
\end{itemize}

\item We claim that $f_{3}$ is injective.

Suppose that $f_{3}(l_{3})=f_{3}(l_{3}^{\prime })$ for some $%
l_{3},l_{3}^{\prime }\in L_{3}.$ Since $\beta _{1}$ is surjective, there
exist $l_{2},l_{2}^{\prime }\in L_{2}$ such that $\beta _{1}(l_{2})=l_{3}$
and $\beta _{1}(l_{2}^{\prime })=l_{3}^{\prime }.$ Then we have%
\begin{equation*}
\begin{array}{rcll}
(f_{3}\circ \beta _{1})(l_{2}) & = & (f_{3}\circ \beta _{1})(l_{2}^{\prime })
&  \\
(\beta _{2}\circ f_{2})(l_{2}) & = & (\beta _{2}\circ f_{2})(l_{2}^{\prime })
&  \\
f_{2}(l_{2})+\alpha _{2}(m_{1}) & = & f_{2}(l_{2}^{\prime })+\alpha
_{2}(m_{1}^{\prime })\text{ (2nd column is exact)} &  \\
(g_{2}\circ \alpha _{2})(m_{1}) & = & (g_{2}\circ \alpha _{2})(m_{1}^{\prime
})\text{ (}g_{2}\circ f_{2}=0\text{)} &  \\
(\alpha _{3}\circ g_{1})(m_{1}) & = & (\alpha _{3}\circ g_{1})(m_{1}^{\prime
}) &  \\
g_{1}(m_{1}) & = & g_{1}(m_{1}^{\prime })\text{ (}\alpha _{3}\text{ is
injective)} &  \\
m_{1}+f_{1}(l_{1}) & = & m_{1}^{\prime }+f_{1}(l_{1}^{\prime })\text{ (1st
row is exact)} &  \\
\alpha _{2}(m_{1})+(\alpha _{2}\circ f_{1})(l_{1}) & = & \alpha
_{2}(m_{1}^{\prime })+(\alpha _{2}\circ f_{1})(l_{1}^{\prime }) &  \\
\alpha _{2}(m_{1})+(f_{2}\circ \alpha _{1})(l_{1}) & = & \alpha
_{2}(m_{1}^{\prime })+(f_{2}\circ \alpha _{1})(l_{1}^{\prime }) &  \\
f_{2}(l_{2})+\alpha _{2}(m_{1})+(f_{2}\circ \alpha _{1})(l_{1}) & = &
f_{2}(l_{2})+\alpha _{2}(m_{1}^{\prime })+(f_{2}\circ \alpha
_{1})(l_{1}^{\prime }) &  \\
f_{2}(l_{2}^{\prime })+\alpha _{2}(m_{1}^{\prime })+(f_{2}\circ \alpha
_{1})(l_{1}) & = & f_{2}(l_{2})+\alpha _{2}(m_{1}^{\prime })+(f_{2}\circ
\alpha _{1})(l_{1}^{\prime }) &  \\
f_{2}(l_{2}^{\prime }+\alpha _{1}(l_{1})) & = & f_{2}(l_{2}+\alpha
_{1}(l_{1}^{\prime }))\text{ (}\alpha _{2}\text{ is cancellative)} &  \\
l_{2}^{\prime }+\alpha _{1}(l_{1}) & = & l_{2}+\alpha _{1}(l_{1}^{\prime })%
\text{ (}f_{2}\text{ is injective)} &  \\
\beta _{1}(l_{2}^{\prime }) & = & \beta _{1}(l_{2})\text{ (}\beta _{1}\circ
\alpha _{1}=0\text{)} &  \\
l_{3}^{\prime } & = & l_{3}.\blacksquare  &
\end{array}%
\newline
\end{equation*}
\end{enumerate}
\end{Beweis}

\begin{proposition}
\label{9}\emph{(The Nine Lemma) }Consider the following commutative diagram
with exact columns and assume that the second row is exact, $\alpha
_{2},f_{2}$ are cancellative and $f_{3},g_{1}$ are $i$-uniform%
\begin{equation*}
\xymatrix{ & 0 \ar@{.>}[d] & 0 \ar[d] & 0 \ar[d] & \\ 0 \ar@{.>}[r] & L_1
\ar[r]^{f_1} \ar[d]_{\alpha_1} & M_1 \ar[r]^{g_1} \ar[d]_{\alpha_2} & N_1
\ar[r] \ar[d]_{\alpha_3} & 0 \\ 0 \ar[r] & L_2 \ar[r]^{f_2} \ar[d]_{\beta_1}
& M_2 \ar[r]^{g_2} \ar[d]_{\beta_2} & N_2 \ar[r] \ar[d]_{\beta_3} & 0 \\ 0
\ar[r] & L_3 \ar[r]^{f_3} \ar[d] & M_3 \ar[r]^{g_3} \ar[d] & N_3
\ar@{-->}[r] \ar@{-->}[d] & 0 \\ & 0 & 0 & 0}
\end{equation*}%
Then the first row is exact if and only if the third row is exact.
\end{proposition}

\begin{proposition}
\label{snake}\emph{(The Snake Lemma) }Consider the following diagram of
semimodules in which the two middle squares are commutative and the two
middle rows are exact. Assume also that the columns are exact (or more
generally that $\alpha _{1},\alpha _{3}$ are $k$-uniform and $\alpha _{2}$
is uniform)%
\begin{equation*}
\xymatrix{ & 0 \ar[d] & 0 \ar[d] & 0 \ar[d] & \\ & {\rm Ker}(\alpha_1)
\ar[d]_{{\rm ker}(\alpha_1)} \ar@{.>}[r]^{f_K} & {\rm Ker}(\alpha_2)
\ar[d]_{{\rm ker}(\alpha_2)} \ar@{.>}[r]^{g_K} & {\rm Ker}(\alpha_3)
\ar[d]_{{\rm ker}(\alpha_3)} \ar@{-->}[dddll]^{\delta} & \\ & L_1
\ar[r]^{f_1} \ar[d]_{\alpha_1} & M_1 \ar[r]^{g_1} \ar[d]_{\alpha_2} & N_1
\ar[r] \ar[d]_{\alpha_3} & 0 \\ 0 \ar[r] & L_2 \ar[r]^{f_2} \ar[d]_{{\rm
coker}(\alpha_1)} & M_2 \ar[r]^{g_2} \ar[d]_{{\rm coker}(\alpha_2)} & N_2
\ar[d]_{{\rm coker}(\alpha_3)} & \\ & {\rm Coker}(\alpha_1)
\ar@{.>}[r]_{f_C} \ar[d] & {\rm Coker}(\alpha_2) \ar@{.>}[r]_{g_C} \ar[d] &
{\rm Coker}(\alpha_3) \ar[d] & \\ & 0 & 0 & 0}
\end{equation*}

\begin{enumerate}
\item There exist unique morphisms $f_{K},g_{K},f_{C}$ and $g_{C}$ which
extend the diagram commutatively.

\item If $f_{1}$ is cancellative, then the first row is exact.

\item If $f_{C}$ is $i$-uniform, then the last row is exact.

\item There exists a $k$-uniform \emph{connecting morphism} $\delta :\mathrm{%
Ker}(\alpha _{3})\longrightarrow \mathrm{Coker}(\alpha _{1})$ such that $%
\mathrm{Ker}(\delta )=\overline{g_{K}(\mathrm{Ker}(\alpha _{2}))}$ and $%
\delta (\mathrm{Ker}(\alpha _{3}))=\mathrm{Ker}(f_{C}).$

\item If $\alpha _{2}$ is cancellative and $g_{K}$ is $i$-uniform, then the
following sequence is exact%
\begin{equation*}
\xymatrix{ & {\rm Ker}(\alpha_2) \ar@{.>}[r]^{g_K} & {\rm Ker}(\alpha_3)
\ar@{-->}[r]^{\delta} & {\rm Coker}(\alpha_1) \ar@{.>}[r]^{f_C} & {\rm
Coker}(\alpha_2) &}
\end{equation*}
\end{enumerate}
\end{proposition}

\begin{Beweis}
\begin{enumerate}
\item The existence and uniqueness of the morphisms $f_{K},g_{K},f_{C}$ and $%
g_{C}$ is guaranteed by the definition of the (co)kernels and the
commutativity of the middle two squares.

\item This follows from Lemma \ref{9-1} applied to the first three rows.

\item This follows from Lemma \ref{9-3} applied to the last three rows.

\item We show first that $\delta $ exists and is well-defined.

\begin{itemize}
\item We define $\delta $ as follows. Let $k_{3}\in \mathrm{Ker}(\alpha
_{3}).$ Choose $m_{1}\in M_{1}$ and $l_{2}\in L_{2}$ such that $%
g_{1}(m_{1})=k_{3}$ and $f_{2}(l_{2})=\alpha _{2}(m_{1});$ notice that this
is possible since $g_{1}$ is surjective and $(g_{2}\circ \alpha
_{2})(m_{1})=(\alpha _{3}\circ g_{1})(m_{1})=\alpha _{3}(k_{3})=0$ whence $%
\alpha _{2}(m_{1})\in \mathrm{Ker}(g_{2})=f_{2}(L_{2}).$ Define $\delta
(k_{3}):=\mathrm{coker}(\alpha _{1})(l_{2})=[l_{2}],$ the coset of $%
L_{2}/\alpha _{1}(L_{1})$ which contains $l_{2}.$

\item $\delta $ is well-defined, i.e. $\delta (k_{3})$ is independent of our
choice of $m_{1}\in M_{1}$ and $l_{2}\in L_{2}$ satisfying the stated
conditions.

Suppose that $g_{1}(m_{1})=k_{3}=g_{1}(m_{1}^{\prime }).$ Since the second
row is exact, there exist $l_{1},l_{1}^{\prime }\in L_{1}$ such that $%
m_{1}+f_{1}(l_{1})=m_{2}+f_{1}(l_{1}^{\prime }).$ It follows that%
\begin{equation*}
\begin{array}{rclc}
\alpha _{2}(m_{1})+(\alpha _{2}\circ f_{1})(l_{1}) & = & \alpha
_{2}(m_{1}^{\prime })+(\alpha _{2}\circ f_{1})(l_{1}^{\prime }) &  \\
f_{2}(l_{2})+(f_{2}\circ \alpha _{1})(l_{1}) & = & f_{2}(l_{2}^{\prime
})+(f_{2}\circ \alpha _{1})(l_{1}^{\prime }) &  \\
f_{2}(l_{2}+\alpha _{1}(l_{1})) & = & f_{2}(l_{2}^{\prime }+\alpha
_{1}(l_{1}^{\prime })) &  \\
l_{2}+\alpha _{1}(l_{1}) & = & l_{2}^{\prime }+\alpha _{1}(l_{1}^{\prime })
& \text{(}f_{2}\text{ is injective)} \\
\lbrack l_{2}] & = & [l_{2}^{\prime }] &
\end{array}%
\end{equation*}%
Thus $l_{2}$ and $l_{2}^{\prime }$ lie in the same coset of $L_{2}/\alpha
_{1}(L_{1}),$ \emph{i.e.} $\delta $ is well-defined.

\item Clearly $\overline{g_{K}(\mathrm{Ker}(\alpha _{2}))}\subseteq \mathrm{%
Ker}(\delta )$ (notice that $f_{2}$ is a semi-monomorphism). We claim that $%
\overline{g_{K}(\mathrm{Ker}(\alpha _{2}))}=\mathrm{Ker}(\delta ).$

Suppose that $k_{3}\in \mathrm{Ker}(\delta )$ for some $k_{3}\in \mathrm{Ker}%
(\alpha _{3}).$ Let $m_{1}\in M_{1}$ be such that $g_{1}(m_{1})=k_{3}$ and
consider $l_{2}\in L_{2}$ such that $f_{2}(l_{2})=\alpha _{2}(m_{1}).$ By
assumption, $[l_{2}]=\delta (k_{3})=0,$ \emph{i.e.} $l_{2}+\alpha
_{1}(l_{1})=\alpha _{1}(l_{1}^{\prime })$ for some $l_{1},l_{1}^{\prime }\in
L_{1}.$Then we have%
\begin{equation*}
\begin{array}{rclc}
f_{2}(l_{2})+(f_{2}\circ \alpha _{1})(l_{1}) & = & (f_{2}\circ \alpha
_{1})(l_{1}^{\prime }) &  \\
\alpha _{2}(m_{1})+\alpha _{2}(f_{1}(l_{1})) & = & \alpha
_{2}(f_{1}(l_{1}^{\prime })) &  \\
m_{1}+f_{1}(l_{1})+k_{2} & = & f_{1}(l_{1}^{\prime })+k_{2}^{\prime } &
\text{(}\alpha _{2}\text{ is }k\text{-uniform)} \\
k_{3}+g_{K}(k_{2}) & = & g_{K}(k_{2}^{\prime }) & \text{(}g_{1}\circ f_{1}=0%
\text{)}%
\end{array}%
\end{equation*}%
Consequently, $\overline{g_{K}(\mathrm{Ker}(\alpha _{2}))}=\mathrm{Ker}%
(\delta ).$

\item Notice that for any $k_{3}\in \mathrm{Ker}(\alpha _{3}),$ we have $%
(f_{C}\circ \delta )(k_{3})=f_{C}([l_{2}])$ where $g_{1}(m_{1})=k_{3}$ and $%
f_{2}(l_{2})=\alpha _{2}(m_{1}).$ It follows that
\begin{equation*}
(f_{C}\circ \delta )(k_{3})=f_{C}(l_{2})=[f_{2}(l_{2})]=[\alpha
_{2}(m_{1})]=[0].
\end{equation*}%
Consequently, $\delta (\mathrm{Ker}(\alpha _{3}))\subseteq \mathrm{Ker}%
(f_{C}).$ We claim that $\delta (\mathrm{Ker}(\alpha _{3}))=\mathrm{Ker}%
(f_{C}).$

Let $[l_{2}]\in \mathrm{Ker}(f_{2}),$ \emph{i.e.} $%
[f_{2}(l_{2})]=f_{C}([l_{2}])=[0],$ for some $l_{2}\in L_{2}.$ Then there
exist $m_{1},m_{1}^{\prime }\in M_{1}$ such that $f_{2}(l_{2})+\alpha
_{2}(m_{1})=\alpha _{2}(m_{1}^{\prime }).$ By assumption, $\alpha _{2}$ is $%
i $-uniform, whence there exists $\mathbf{m}_{1}\in M_{1}$ such that $\alpha
_{2}(\mathbf{m}_{1})=f_{2}(l_{2}).$ It follows that $(\alpha _{3}\circ
g_{1})(\mathbf{m}_{1})=(g_{2}\circ \alpha _{2})(\mathbf{m}_{1})=(g_{2}\circ
f_{2})(l_{2})=0.$ So, $g_{1}(\mathbf{m}_{1})\in \mathrm{Ker}(\alpha _{3})$
and $\delta (g_{1}(\mathbf{m}_{1}))=[l_{2}].$ Consequently, $\mathrm{Ker}%
(f_{C})=\delta (\mathrm{Ker}(\alpha _{3})).$

\item We claim that $\delta $ is $k$-uniform.

Suppose that $\delta (k_{3})=\delta (k_{3}^{\prime })$ for some $%
k_{3},k_{3}^{\prime }\in \mathrm{Ker}(\alpha _{3}).$ Let $%
m_{1},m_{1}^{\prime }\in M_{1}$ and $l_{2},l_{2}^{\prime }\in L_{2}$ be such
that $g_{1}(m_{1})=k_{3},$ $g_{1}(m_{1}^{\prime })=k_{3}^{\prime },$ $\alpha
_{2}(m_{1})=f_{2}(l_{2})$ and $\alpha _{2}(m_{1}^{\prime
})=f_{2}(l_{2}^{\prime }).$ By assumption, $[l_{2}]=[l_{2}^{\prime }],$
\emph{i.e. }$l_{2}+\alpha _{1}(l_{1})=l_{2}+\alpha _{1}(l_{1}^{\prime })$
for some $l_{1},l_{1}^{\prime }\in L_{1}.$ Notice that%
\begin{equation*}
\begin{array}{rclc}
f_{2}(l_{2})+(f_{2}\circ \alpha _{1})(l_{1}) & = & f_{2}(l_{2}^{\prime
})+(f_{2}\circ \alpha _{1})(l_{1}^{\prime }) &  \\
\alpha _{2}(m_{1})+(\alpha _{2}\circ f_{1})(l_{1}) & = & \alpha
_{2}(m_{2}^{\prime })+(\alpha _{2}\circ f_{1})(l_{1}^{\prime }) &  \\
m_{1}+f_{1}(l_{1})+k_{2} & = & m_{1}^{\prime }+f_{1}(l_{1}^{\prime
})+k_{2}^{\prime } & \text{(}\alpha _{2}\text{ is }k\text{-uniform)} \\
g_{1}(m_{1})+g_{K}(k_{2}) & = & g_{1}(m_{1}^{\prime })+g_{K}(k_{2}^{\prime })
& \text{(}g_{1}\circ f_{1}=0\text{)} \\
k_{3}+g_{K}(m_{1}) & = & k_{3}^{\prime }+g_{K}(m_{1}) &
\end{array}%
\end{equation*}%
\qquad Since $g_{K}(\mathrm{Ker}(\alpha _{2}))\subseteq \mathrm{Ker}(\delta
) $ we conclude that $\delta $ is $k$-uniform.
\end{itemize}

\item If $g_{K}$ is $i$-uniform, then we have $\mathrm{Ker}(\delta )=%
\overline{g_{K}(\mathrm{Ker}(\alpha _{2}))}=g_{K}(\mathrm{Ker}(\alpha _{2}))$
and it remains only to prove that $f_{C}$ is $k$-uniform.

Suppose that $f_{C}[l_{2}]=f_{C}[l_{2}^{\prime }]$ for some $%
l_{2},l_{2}^{\prime }\in L_{2}.$ Then there exist $m_{1},m_{1}^{\prime }\in
M_{1}$ such that $f_{2}(l_{2})+\alpha _{2}(m_{1})=f_{2}(l_{2}^{\prime
})+\alpha _{2}(m_{1}^{\prime }).$ It follows that%
\begin{equation*}
\begin{array}{rclc}
(g_{2}\circ \alpha _{2})(m_{1}) & = & (g_{2}\circ \alpha _{2})(m_{1}^{\prime
})\text{ (}g_{2}\circ f_{2}=0\text{)} &  \\
(\alpha _{3}\circ g_{1})(m_{1}) & = & (\alpha _{3}\circ g_{1})(m_{1}^{\prime
}) &  \\
g_{1}(m_{1})+k_{3} & = & g_{1}(m_{1}^{\prime })+k_{3}^{\prime }\text{ (}%
\alpha _{3}\text{ is }k\text{-uniform)} &  \\
g_{1}(m_{1}+\mathbf{m}_{1}) & = & g_{1}(m_{1}^{\prime }+\mathbf{m}%
_{1}^{\prime })\text{ (}g_{1}\text{ is surjective)} &  \\
m_{1}+\mathbf{m}_{1}+f_{1}(\mathbf{l}_{1}) & = & m_{1}^{\prime }+\mathbf{m}%
_{1}+f_{1}(\mathbf{l}_{1}^{\prime })\text{ (2nd row is exact)} &  \\
\alpha _{2}(m_{1})+\alpha _{2}(\mathbf{m}_{1})+(\alpha _{2}\circ f_{1})(%
\mathbf{l}_{1}) & = & \alpha _{2}(m_{1}^{\prime })+\alpha _{2}(\mathbf{m}%
_{1}^{\prime })+(\alpha _{2}\circ f_{1})(\mathbf{l}_{1}^{\prime }) &  \\
f_{2}(l_{2}^{\prime })+\alpha _{2}(m_{1})+\alpha _{2}(\mathbf{m}%
_{1})+(f_{2}\circ \alpha _{1})(\mathbf{l}_{1}) & = & [f_{2}(l_{2}^{\prime
})+\alpha _{2}(m_{1}^{\prime })]+\alpha _{2}(\mathbf{m}_{1}^{\prime
})+(f_{2}\circ \alpha _{1})(\mathbf{l}_{1}^{\prime }) &  \\
f_{2}(l_{2}^{\prime })+\alpha _{2}(m_{1})+\alpha _{2}(\mathbf{m}%
_{1})+(f_{2}\circ \alpha _{1})(\mathbf{l}_{1}) & = & f_{2}(l_{2})+\alpha
_{2}(m_{1})+\alpha _{2}(\mathbf{m}_{1}^{\prime })+(f_{2}\circ \alpha _{1})(%
\mathbf{l}_{1}^{\prime }) &  \\
f_{2}(l_{2}^{\prime })+\alpha _{2}(\mathbf{m}_{1})+(f_{2}\circ \alpha _{1})(%
\mathbf{l}_{1}) & = & f_{2}(l_{2})+\alpha _{2}(\mathbf{m}_{1}^{\prime
})+(f_{2}\circ \alpha _{1})(\mathbf{l}_{1}^{\prime })\text{ (}\alpha _{2}%
\text{ is cancellative)} &  \\
f_{2}(l_{2}^{\prime }+\mathbf{l}_{2}+\alpha _{1}(\mathbf{l}_{1})) & = &
f_{2}(l_{2}+\mathbf{l}_{2}^{\prime }+\alpha _{1}(\mathbf{l}_{1}^{\prime }))%
\text{ (third row is exact)} &  \\
l_{2}^{\prime }+\mathbf{l}_{2}+\alpha _{1}(\mathbf{l}_{1}) & = & l_{2}+%
\mathbf{l}_{2}^{\prime }+\alpha _{1}(\mathbf{l}_{1}^{\prime })\text{ (}f_{2}%
\text{ is injective)} &  \\
\lbrack l_{2}^{\prime }]+[\mathbf{l}_{2}] & = & [l_{2}]+[\mathbf{l}%
_{2}^{\prime }] &  \\
\lbrack l_{2}^{\prime }]+\delta (k_{3}) & = & [l_{2}]+\delta (k_{3}^{\prime
}) &
\end{array}%
\end{equation*}%
Since $\delta (\mathrm{Ker}(\alpha _{3}))\subseteq \mathrm{Ker}(f_{C}),$ we
conclude that $f_{C}$ is $k$-uniform.$\blacksquare $
\end{enumerate}
\end{Beweis}

\textbf{Acknowledgments.} The author thanks all mathematicians who clarified
to him some issues related to the nature of the categories of semimodules
and exact sequences or sent to him related manuscripts especially F. Linton,
G. Janelidze, Y. Katsov, A. Patchkoria, H. Porst and R. Wisbauer.

\end{document}